\documentclass[12pt,reqno]{amsart} 

\usepackage{float} 

\usepackage{epsfig,amssymb,amsmath,version}
\usepackage{amssymb,version,graphicx,fancybox,mathrsfs,labelfig}

\usepackage{url,hyperref}
\usepackage{subfigure}
\usepackage{color}
\usepackage{stmaryrd}
\usepackage{multirow}
\usepackage{booktabs,siunitx}
\usepackage{booktabs}
\usepackage{multicol,epstopdf}
\usepackage{xypic}
\usepackage[]{graphicx}
\usepackage[all]{xy}
\usepackage{tikz,ifthen}
\usetikzlibrary{positioning, math, decorations.markings, arrows.meta}
\usetikzlibrary{calc,shapes,cd}
\usetikzlibrary{knots} 
\usepgfmodule{decorations}
\allowdisplaybreaks
\tikzset{knotarrow/.pic={ \draw[edge, <-] (0,0) -- +(-.001,0);}}

\tikzset{edge/.style={line width=0.8}}
\tikzset{wall/.style={very thick}}

\tikzset{->-/.style n args={2}{decoration={markings, mark=at position #1 with {\arrow{#2}}}, postaction={decorate}}} 

\ExplSyntaxOn
\cs_new_eq:NN \ifstreqF \str_if_eq:nnF
\cs_new_eq:NN \ifstreqTF \str_if_eq:nnTF
\ExplSyntaxOff

\tikzset{-o-/.code 2 args={\ifstreqF{#2}{} 
{\ifstreqTF{#2}{>}
   {\pgfkeysalso{decoration={markings,mark=at position #1 with {\arrow[scale=0.8]{#2}}}
                    ,postaction={decorate}}
    }
   {\ifstreqTF{#2}{<}
       {\pgfkeysalso{decoration={markings,mark=at position #1 with {\arrow[scale=0.8]{#2}}}
                    ,postaction={decorate}}
        }
       {\pgfkeysalso{decoration={markings,
                    mark=at position #1 with
                    {\draw[black, fill={#2}] circle[radius=2pt];}}
                    ,postaction={decorate}}
        }
     }
  }}}

\allowdisplaybreaks[4]

\newtheorem{theorem}{Theorem}[section]
\newtheorem{lemma}[theorem]{Lemma}

\newtheorem{corollary}[theorem]{Corollary}
\newtheorem{proposition}[theorem]{Proposition}
\newtheorem{remark}[theorem]{Remark}

\def \cred {\color{red}}

\definecolor{ligreen}{rgb}{0.0, 0.3, 0.0}
\def \cg {\color{ligreen}}

\definecolor{darkblue}{rgb}{0.0, 0.0, 0.55}

\definecolor{anti-flashwhite}{rgb}{0.55, 0.57, 0.68}
\def \cF {\mathcal{F}}
\def \Sq {S_q(M)}
\def \Sqp {S_q(\partial M)}

\def \Sz {S_{\zeta}(M)}
\def \Szp {S_{\zeta}(\partial M)}
\def \SqN {S_q(M)^{(N)}}

\def \Si {S_q(\Sigma)}
\def \Siz {S_{\zeta}(\Sigma)}
\def \SiN {S_q(\Sigma)^{(N)}}

\def \pM {\partial M}

\def \define {\triangleq}
\def \Ls {L_{\sharp}}

\def \Tu {\mathcal{T}^{\mu}_{\tau}}
\def \Tv {\mathcal{T}^{\nu}_{\tau}}

\def \Zu {\mathcal{Z}^{\mu}_{\tau}}
\def \Zv {\mathcal{Z}^{\nu}_{\tau}}

\def \Tq {Tr_{\tau}^{q}}
\def \Tz {Tr_{\tau}^{\zeta}}

\def \bZ{\mathbb{Z}}

\newcommand{\bp}{\begin{proposition}}
\newcommand{\ep}{\end{proposition}}
\newcommand{\bpr}{\begin{proof}}
\newcommand{\epr}{\end{proof}}

\newcommand{\bt}{\begin{theorem}}
\newcommand{\et}{\end{theorem}}

\newcommand{\bl}{\begin{lemma}}
\newcommand{\el}{\end{lemma}}

\newcommand{\bcr}{\begin{corollary}}
\newcommand{\ecr}{\end{corollary}}

\newcommand{\be}{\begin{equation}}
\newcommand{\ee}{\end{equation}}

\newcommand{\bes}{\begin{equation*}}
\newcommand{\ees}{\end{equation*}}

\newcommand{\ba}{\begin{align}}
\newcommand{\ea}{\end{align}}

\newcommand{\bas}{\begin{align*}}
\newcommand{\eas}{\end{align*}}

\graphicspath{{./Figures/} } 

\begin{document}
\bibliographystyle{plain}

\title{Finiteness and dimension of  stated skein modules over  Frobenius}
\author{Zhihao Wang}


 \maketitle


\begin{abstract}
When the quantum parameter $q^{1/2}$ is a root of unity of odd order,
the stated skein module $S_{q^{1/2}}(M,\mathcal{N})$ has an $S_{1}(M,\mathcal{N})$-module structure, where $(M,\mathcal{N})$ is a marked three manifold.
We prove
$S_{q^{1/2}}(M,\mathcal{N})$ is a finitely generated  $S_{1}(M,\mathcal{N})$-module when $M$ is compact, which furthermore indicates the reduced stated skein module for the compact marked three manifold is finite dimensional. We also give an upper bound for the dimension of $S_{q^{1/2}}(M,\mathcal{N})$ over  $S_{1}(M,\mathcal{N})$ when $M$ is compact.

For a punctured bordered surface $\Sigma$, we use $S_{q^{1/2}}(\Sigma)^{(N)}$ to denote the image of the Frobenius map when $q^{1/2}$ is a root of unity of odd order $N$. Then $S_{q^{1/2}}(\Sigma)^{(N)}$ lives in the center of the stated skein algebra $S_{q^{1/2}}(\Sigma)$. Let $\widetilde{S_{q^{1/2}}(\Sigma)^{(N)}}$ be the field of fractions of $S_{q^{1/2}}(\Sigma)^{(N)}$, and $\widetilde{S_{q^{1/2}}(\Sigma)}$ be $S_{q^{1/2}}(\Sigma)\otimes_{S_{q^{1/2}}(\Sigma)^{(N)}} \widetilde{S_{q^{1/2}}(\Sigma)^{(N)}}$. Then we show the dimension of $\widetilde{S_{q^{1/2}}(\Sigma)}$ over 
$\widetilde{S_{q^{1/2}}(\Sigma)^{(N)}}$ is $N^{3r(\Sigma)}$ where $r(\Sigma)$ equals to  the number of boundary components of $\Sigma$ minus the Euler characteristic of $\Sigma$.
\end{abstract}


\def \SMQ {S_{q^{1/2}}(M,\mathcal{N})}
\def \SMQP {S_{q^{1/2}}(M^{'},\mathcal{N}^{'})}

\def \CSM {S_1(M,\mathcal{N})}

\def \SMN {S_{q^{1/2}}(M,\mathcal{N})^{(N)}}

\def \S {S_{q^{1/2}}}

\def \MN {(M,\mathcal{N})}

\section{Introduction}


In this paper, our ground ring is the complex field $\mathbb{C}$ with a distinguished nonzero element $q^{\frac{1}{2}}$. 

All the three manifolds and surfaces mentioned in this paper are assumed to be oriented.


 A {\bf marked three manifold} is  a pair $(M,\mathcal{N})$, 
where $M$ is a three manifold, and $\mathcal{N}$ is a one dimensional submanifold of $\partial M$ consisting of oriented open intervals such that there is no intersection between the closure of any two open intervals. Note that it is possible that $\mathcal{N} = \emptyset$. We say $\MN$ is compact if $M$ is compact.


The
stated skein module is a generalization for the Kauffman bracket skein module \cite{przytycki2006skein,turaev1988conway}. It was first defined by Bonahon and Wong
to establish the famous quantum trace map \cite{bonahon2011quantum}, then it was refined by L{\^e} \cite{le2018triangular}. The precise definition of the stated skein module for a marked three manifold $\MN$,   denoted as $\S\MN$, came from \cite{bloomquist2020chebyshev}. Compared with the skein theory, the stated case has more structures, such as comodule structure over $O_q(SL_2)$, splitting map, etc \cite{costantino2022stated1,costantino2022stated}.

When $q^{\frac{1}{2}}$ is a root of unity of order $N$ with $N$ being odd,
 there is an important linear map, called Frobenius map, $\cF:\CSM\rightarrow \SMQ$ \cite{bloomquist2020chebyshev,bonahon2016representations,wang2023stated}.
This map is helpful to understand the center of the stated skein algebra, and the representation theory of the stated skein algebra \cite{bloomquist2020chebyshev,bonahon2016representations,bonahon2017representations,frohman2019unicity,korinman2021unicity,korinman2019classical,wang2023stated}.
This map is also the key to define the reduced skein module \cite{detcherry2023kauffman,talk}, which was used to calculate the dimension of the skein module for some closed three manifolds \cite{detcherry2023kauffman}.

\def \SMn {S_{q^{1/2}}(M,\mathcal{N})^{(N)}}

$S_1(M,\mathcal{N})$ has a commutative algebra structure given by taking disjoint union of stated tangles in $\MN$.
Let $\SMn=\text{Im}\mathcal{F}$. 
We show $\SMn$ also has a commutative algebra structure, which makes 
$\cF:\CSM\rightarrow \SMn$ a surjective algebra homomorphism.

\def \SMN {S_{1}(M,\mathcal{N})}
As explained later, 
$\SMQ$ has an $\SMN$-module structure. It was proved that $\SMQ$ is finitely generated over $\SMN$ when
$\MN$ is the thickening of a pb surface (see definition for the pb surface in subsection \ref{2.1}), which actually is the crucial step to prove the Unicity Theorem for stated skein algebras \cite{abdiel2017localized,frohman2019unicity,korinman2021unicity,wang2023stated}. 
In this paper, we will prove the following more general result.

\begin{theorem}\label{1.1}
When $q^{\frac{1}{2}}$ is a root of unity of order $N$ with $N$ being odd,
 $\SMQ$ is finitely generated over $\SMN$ for any compact marked three manifold $\MN$. 
\end{theorem}

We call the property that $\SMQ$ is finitely generated over $\SMN$ as  finiteness of  stated skein modules over  Frobenius.

 Theorem \ref{1.1} was conjectured in \cite{wang2023stated} for general stated $SL_n$-skein modules.
Although we will only work with the $SL_2$ case for convenience, some
 our techniques used in this paper can be easily generalized to the general $SL_n$ case. We will state  corresponding general results for the $SL_n$ case without giving detailed proof.
Theorem \ref{1.1}  implies the following  Theorem about the reduced stated skein module, defined in subsection \ref{sub_reduced}.
\bt
The
reduced stated skein module of a compact marked three manifold is finite dimensional. 
\et

\def \cF {\mathcal{F}}
For a three  manifold $M$, we will use $\Sq$ to denote the Kauffman bracket skein module of $M$.
Let $q^2$ be a root of unity of  order $N$ with $N$ being odd, and $M$ be a three manifold.
Then the Frobenius map $\cF: S_{q^N}(M)\rightarrow S_q(M)$ is defined (note  that $q^N = \pm1$) \cite{bloomquist2020chebyshev,bonahon2016representations}.

Since $q^N =\pm 1$, $S_{q^N}(M)$ has a commutative algebra structure given by taking disjoint union of framed links.
 Similarly we have 
 $S_q(M)$ has an 
$S_{q^N}(M)$-module structure. So we have a stronger version for Theorem \ref{1.1} when $\mathcal{N}=\emptyset$. That is,
$S_q(M)$ is finitely generated over $S_{q^N}(M)$ when $M$ is compact and $q^2$ is a root of unity of  order $N$ with $N$ being odd \cite{talk}.

For every three manifold $M$, we know the skein algebra $S_q(\partial M)$ acts on $\Sq$.
Using the stronger version for Theorem \ref{1.1} when $\mathcal{N}=\emptyset$, 
we find out finiteness of $S_{q^N}(M)$ over  $S_{q^N}(\partial M)$ implies finiteness of $\Sq$ over  $S_q(\partial M)$ when $M$ is compact and $q^2$ is a root of unity of  order $N$ with $N$ being odd.

\bt\label{thm1.3}
Let $M$ be a compact   three manifold. Suppose $q^2$ is a root of unity of order $N$ with $N$ being odd, and
 $S_{q^N}(M)$ is finitely generated over $S_{q^N}(\partial M)$.
 Then $\Sq$ is finitely generated over $\Sqp$. 
\et

For a compact marked three manifold $\MN$,
we  also give an upper bound for the dimension of $\SMQ$ over $\SMN$, which is defined to be the minimal number of generators of $\SMQ$ over $\SMN$. This furthermore gives an upper bound for the dimension of the reduced stated skein module of $\MN$.
\bt\label{dimensi}
Suppose
$q^{\frac{1}{2}}$ is a root of unity of order $N$ with $N$ being odd.
Let $(M,\mathcal{N})$ be a connected compact marked three manifold. 
Suppose the Heegaard genus of $M$ is $g$, and $\mathcal{N}$ has $k$ components.

(a) When $k=0$, we have the dimension of $S_{q^{1/2}}(M,\emptyset)$ over $S_1(M,\emptyset)$ is not more than 
 $ N^{2^g-1}$.

(b) When $k>0$, we have 
  the dimension of $\SMQ$ over $\SMN$ is not more than 
$$ (2N^3 - \frac{N(N+1)(2N+1)}{6})^{2g+k-1}.$$
\et

\def \q {q^{\frac{1}{2}}}

\def \SqSN {S_{q^{1/2}}(\Sigma)^{(N)}}
\def \fSiN {\widetilde{S_{q^{1/2}}(\Sigma)^{(N)}}}
\def \fSi {\widetilde{\S(\Sigma)}}

Let $\Sigma$ be a pb surface,  we know that $\S(\Sigma)$ is a domain. Then when $\q$ is a root of unity of order $N$ with $N$ being odd,
 $\S(\Sigma)^{(N)} =\text{Im}\mathcal{F}$ is a commutative domain. 
We use $\fSiN$ to denote the field of fractions of $\S(\Sigma)^{(N)}$, and use 
$\fSi$ to denote $\S(\Sigma)\otimes_{\S(\Sigma)^{(N)}}\fSiN$. Then $\fSi$ is a finite dimensional vector space over the field
$\fSiN$. We define 
$r(\Sigma) = -\chi(\Sigma)+\sharp\partial\Sigma$, where $\chi(\Sigma)$ is the Euler characteristic of $\Sigma$ and $\sharp\partial\Sigma$ is the number of boundary components of $\Sigma$.
\bt\label{thm1.5}
Suppose
$q^{\frac{1}{2}}$ is a root of unity of order $N$ with $N$ being odd.
Let $\Sigma$ be a connected  pb surface, we require $\chi(\Sigma)$ is negative if $\partial \Sigma = \emptyset$.  
We have $$\text{dim}_{\fSiN}\fSi= N^{3 r(\Sigma)}.$$
\et

When $\partial \Sigma =\emptyset$, Theorem \ref{thm1.5} was already proved in \cite{frohman2018structure}.   Theorem 6.1 in \cite{frohman2021dimension} also  indicates the case when $\partial \Sigma=\emptyset$ in Theorem \ref{thm1.5}.

Thang L{\^e} and Tao Yu claimed the center  and the PI-dimension of the stated skein algebra in
\cite{le2021statedsss} (they allow $N$ to be even). When the center is the image of the Frobenius map, the PI-dimension should be the square root of $\text{dim}_{\fSiN}\fSi$. It is a trivial check that  Theorem \ref{thm1.5} coincides with  Theorem 7.3 in \cite{le2021statedsss} when the center of the stated skein algebra is the image of the Frobenius map.


%

%

%


{\bf Plan of this paper.} In section 2, we introduce some important definitions, such as stated
skein modules (algebras), Chebyshev polynomials, and the Frobenius map. In section 3,  we will show the connection between Theorem \ref{1.1} and the finiteness of the skein module of the three manifold $M$ over the skein algebra of $\partial M$, which indicates Theorem \ref{thm1.3}. In section 4, we prove Theorem \ref{1.1}. 
In section 5, we prove the center of the skein algebra is freely generated by $N^p$ elements over the image of the Frobenius map using the quantum trace map, where $p$ is the number of punctures. In section 6, we  prove Theorems \ref{dimensi} and \ref{thm1.5}.

{\bf
Acknowledgements}: The author would like to thank Thang L{\^e} for helpful comments. The research is supported by NTU  research scholarship, and PhD  scholarship from University of Groningen.

\section{Preliminary}

We use $\mathbb{N}$ to denote the set of nonnegative integers. and use $\mathbb{Z}$ to denote the set of integers. We set $\mathbb{N}^{*} = \mathbb{N}-\{0\}, \mathbb{Z}^{*} = \mathbb{Z}-\{0\}, \mathbb{C}^{*} = \mathbb{C}-\{0\}$.

Let $A$ be an algebra,  we use $C(A)$ to denote the center of $A$. Suppose $B$ is a finitely generated $A$-module. We use dim${}_A B$ to denote the lower bound of the number of generators of $B$ over $A$. 
We call an algebra $A$ a domain, if $xy = 0$ for  two elments $x,y\in A$, we have $x=0$ or $y=0$.

\bl\label{free}
Suppose $A$ is a commutative domain, and $B$ is a module over $A$ freely generated by $n$ elements. Then dim${}_A B=n$.
\el
\bpr
Clearly we have dim${}_A B\leq n$. It suffices to show dim${}_A B\geq n$. Let $S = A-\{0\}$. 

We use 
$\tilde{A}$ to denote the localization of $A$ over the multiplicitive set $S$, that is, $\tilde{A}$ is the field of fractions of $A$. Then 
$B\otimes _{A}\tilde{A}$ is a vertor space over $\tilde{A}$. From the assumption, we have $B\simeq A^{\oplus n}$ as $A$-modules. Then we have 
\be\label{iso}
B\otimes _{A}\tilde{A}\simeq ( A^{\oplus n})\otimes _{A}\tilde{A}
= (A\otimes _{A}\tilde{A})^{\oplus n}\simeq \tilde{A}^{\oplus n}
\ee
where all the isomorphisms in equation \eqref{iso} are isomorphisms between $\tilde{A}$-vector spaces.
Thus dim$_{\tilde{A}} (B\otimes _{A}\tilde{A}) = n$. Let $\{b_1,\dots, b_k\}$ be any finite generating set for $B$ over $A$. Then we have $\{b_1\otimes _{A} 1,\dots, b_k\otimes _{A} 1\}$ linearly spans $B\otimes _{A}\tilde{A}$, thus we have dim$_{\tilde{A}}( B\otimes _{A}\tilde{A}) = n\leq k$. Then we get $n\leq$dim${}_A B$.
\epr

\def \SMQ {S_{q^{1/2}}(M,\mathcal{N})}
\def \SMQP {S_{q^{1/2}}(M^{'},\mathcal{N}^{'})}

\subsection{Stated skein modules for marked three manifolds and stated skein algebras for punctured bordered surfaces and marked surfaces}\label{2.1}
Let $(M,\mathcal{N})$ be a marked three manifold. An {\bf $(M,\mathcal{N})$-tangle} $\alpha$ is a properly embedded one dimensional submanifold of $M$, equipped with a framing such that $\partial \alpha\subset \mathcal{N}$ and at the each point in $\partial \alpha$ the framing is given by the velocity vector of $\mathcal{N}$. Two $(M,\mathcal{N})$-tangles are isotopic if they are isotopic in the class of $(M,\mathcal{N})$-tangles.

A {\bf stated $(M,\mathcal{N})$-tangle} is an $(M,\mathcal{N})$-tangle $\alpha$ equipped with a map $s:\partial\alpha\rightarrow \{-,+\}$.
We will consider stated $(M,\mathcal{N})$-tangles up to isotopy. 

The stated skein module $\SMQ$ is 
defined to be the vector space over $\mathbb{C}$ with the set of all isotopy classes of stated $(M,\mathcal{N})$-tangles as the basis subject to the following relations:
\begin{equation}\label{cross}
\raisebox{-.20in}{
\begin{tikzpicture}
\filldraw[draw=white,fill=gray!20] (-0,-0.2) rectangle (1, 1.2);
\draw [line width =1pt](0.6,0.6)--(1,1);
\draw [line width =1pt](0.6,0.4)--(1,0);
\draw[line width =1pt] (0,0)--(0.4,0.4);
\draw[line width =1pt] (0,1)--(0.4,0.6);
\draw[line width =1pt] (0.6,0.6)--(0.4,0.4);
\end{tikzpicture}
}=
q
\raisebox{-.20in}{
\begin{tikzpicture}
\filldraw[draw=white,fill=gray!20] (-0,-0.2) rectangle (1, 1.2);
\draw [line width =1pt](0.6,0.6)--(1,1);
\draw [line width =1pt](0.6,0.4)--(1,0);
\draw[line width =1pt] (0,0)--(0.4,0.4);
\draw[line width =1pt] (0,1)--(0.4,0.6);
\draw[line width =1pt] (0.6,0.62)--(0.6,0.38);
\draw[line width =1pt] (0.4,0.38)--(0.4,0.62);
\end{tikzpicture}
}
+
 q^{-1}
\raisebox{-.20in}{
\begin{tikzpicture}
\filldraw[draw=white,fill=gray!20] (-0,-0.2) rectangle (1, 1.2);
\draw [line width =1pt](0.6,0.6)--(1,1);
\draw [line width =1pt](0.6,0.4)--(1,0);
\draw[line width =1pt] (0,0)--(0.4,0.4);
\draw[line width =1pt] (0,1)--(0.4,0.6);
\draw[line width =1pt] (0.62,0.6)--(0.38,0.6);
\draw[line width =1pt] (0.62,0.4)--(0.38,0.4);
\end{tikzpicture}
} 
\end{equation}
\begin{equation}\label{unknot}
\raisebox{-.15in}{
\begin{tikzpicture}
\filldraw[draw=white,fill=gray!20] (-0,-0) rectangle (1, 1);
\draw [line width =1pt] (0.5,0.5) circle (0.3);
\end{tikzpicture}
}=-(q^2+q^{-2})
\raisebox{-.15in}{
\begin{tikzpicture}
\filldraw[draw=white,fill=gray!20] (-0,-0) rectangle (1, 1);
\end{tikzpicture}
}
\end{equation}
\begin{equation}\label{arc}
\raisebox{-.26in}{
\begin{tikzpicture}
\filldraw[draw=white,fill=gray!20] (-0,-0) rectangle (1, 1);
\draw [line width =1pt]  (0.5 ,0) arc (-90:225:0.3 and 0.35);
\filldraw[draw=black,fill=black] (0.5,-0) circle (0.09);
\node at (0.3,-0.15) {$-$};
\node at (0.7,-0.15) {\small $+$};
\end{tikzpicture}
}=q^{-\frac{1}{2}}
\raisebox{-.15in}{
\begin{tikzpicture}
\filldraw[draw=white,fill=gray!20] (-0,-0) rectangle (1, 1);
\filldraw[draw=black,fill=black] (0.5,-0) circle (0.09);
\end{tikzpicture}
},\;
\raisebox{-.26in}{
\begin{tikzpicture}
\filldraw[draw=white,fill=gray!20] (-0,-0) rectangle (1, 1);
\draw [line width =1pt]  (0.5 ,0) arc (-90:225:0.3 and 0.35);
\filldraw[draw=black,fill=black] (0.5,-0) circle (0.09);
\node at (0.3,-0.15) {$+$};
\node at (0.7,-0.15) {\small $+$};
\end{tikzpicture}
}=
\raisebox{-.26in}{
\begin{tikzpicture}
\filldraw[draw=white,fill=gray!20] (-0,-0) rectangle (1, 1);
\draw [line width =1pt]  (0.5 ,0) arc (-90:225:0.3 and 0.35);
\filldraw[draw=black,fill=black] (0.5,-0) circle (0.09);
\node at (0.3,-0.15) {$-$};
\node at (0.7,-0.15) {\small $-$};
\end{tikzpicture}
} =0
\end{equation}
\begin{equation}\label{hight}
\raisebox{-.26in}{
\begin{tikzpicture}
\filldraw[draw=white,fill=gray!20] (-0.2,-0) rectangle (1.2, 1);
\draw [line width =1pt](0,1)--(0.4,0.2);
\draw [line width =1pt](1,1)--(0.6,0.2);
\draw [line width =1pt](0.5,0)--(0.6,0.2);
\filldraw[draw=black,fill=black] (0.5,-0) circle (0.09);
\node at (0.3,-0.15) {$+$};
\node at (0.7,-0.15) {\small $-$};
\end{tikzpicture}
}=q^{2}
\raisebox{-.26in}{
\begin{tikzpicture}
\filldraw[draw=white,fill=gray!20] (-0.2,-0) rectangle (1.2, 1);
\draw [line width =1pt](0,1)--(0.4,0.2);
\draw [line width =1pt](1,1)--(0.6,0.2);
\draw [line width =1pt](0.5,0)--(0.6,0.2);
\filldraw[draw=black,fill=black] (0.5,-0) circle (0.09);
\node at (0.3,-0.15) {$-$};
\node at (0.7,-0.15) {\small $+$};
\end{tikzpicture}
} + q^{-\frac{1}{2}}
\raisebox{-.16in}{
\begin{tikzpicture}
\filldraw[draw=white,fill=gray!20] (-0.2,-0) rectangle (1.2, 1);
\draw [line width =1pt](0,1)--(0.4,0.2);
\draw [line width =1pt](1,1)--(0.6,0.2);
\draw [line width =1pt](0.38,0.2)--(0.62,0.2);
\filldraw[draw=black,fill=black] (0.5,-0) circle (0.09);
\end{tikzpicture}
}
\end{equation}
where the black dot represents a component of $\mathcal{N}$, which is perpendicular to this page and points to readers, and the black lines are parts of the stated $(M,\mathcal{N})$-tangles with framing pointing to readers, the gray square is the projection of an embedded cube in $M$, all the stated tangles in a same equation are identical outside the cube.
 For detailed explanation for the above relations and stated skein modules, please refer to \cite{bloomquist2020chebyshev, le2018triangular}.

When $\mathcal{N}$ is empty, $\SMQ$ is just the (Kauffman bracket) skein module. In this case we denote $S_{q^{1/2}}(M,\emptyset)$
as $S_q(M)$, and call $S_q(M)$ the skein module of the three manifold $M$ (the reason why we replace $\q$ with $q$ for skein module is because relations \eqref{cross} and \eqref{unknot} only involve $q$).

Let $(M_1,\mathcal{N}_1), (M_2, \mathcal{N}_2)$ be two marked three manifolds, and 
$f:M_1\rightarrow M_2$ be an orientation preversing embedding. If furthermore 
$f(\mathcal{N}_1)\subset \mathcal{N}_2$ and preserves the orientations of the marking, we say $f$ is an embedding from 
$(M_1,\mathcal{N}_1)$ to $ (M_2, \mathcal{N}_2)$. Obviously $f$ induces a linear map
$f_{\sharp} : \S (M_1,\mathcal{N}_1)\rightarrow \S(M_2,\mathcal{N}_2)$. We say $f$ is an isomorphism if there exists another embedding
$g:(M_2, \mathcal{N}_2)\rightarrow (M_1,\mathcal{N}_1)$ such that $g$ and $f$ are inverse to each other. We also say $(M_1,\mathcal{N}_1)$ is isomorphic to $(M_2, \mathcal{N}_2)$.



A {\bf punctured bordered surface} 
is obtained from a compact surface by removing finite points such that every boundary component of the punctured bordered surface is isomorphic to the open interval $(0,1)$. These removed points are called {\bf punctures}, and the punctures are called {\bf boundary punctures} if they are on the boundary, otherwise they are called {\bf interior punctures}.
 For simplicity we will call punctured bordered surface a {\bf pb surface}. When $\partial\Sigma$ is empty, we say $\Sigma$ is a {\bf closed pb  surface}. A {\bf totally closed pb  surface} is a closed pb surface without punctures.


A {\bf triangulation} of  a pb surface with at least one puncture is a collection of non-isotopic simple arcs connecting punctures such that no arc bounds a disk and these arcs cut the surface into ideal triangles (an ideal triangle is obtained from a disk by removing three points on the boundary).

Let $\Sigma$ be a pb surface. For any component $e$ of $\partial \Sigma$, we select a point $x_e\in e$. Set $M = \Sigma\times[0,1], \mathcal{N} = \cup_{e}( \{x_e\}\times (0,1))$
where the union takes over all components of $\partial\Sigma$, then $(M,\mathcal{N})$ is a marked three manifold, and we say $(M,\mathcal{N})$ is the thickening of $\Sigma$. And we define $S_{q^{1/2}}(\Sigma)$ to be $\SMQ$.
For any two stated $(M,\mathcal{N})$-tangles
$\alpha,\beta$, the product $\alpha\beta$ is obtained by stacking $\alpha$ above $\beta$. This makes 
$S_{q^{1/2}}(\Sigma)$ into an algebra, called the stated skein algebra of $\Sigma$. 


When $\Sigma$ is closed, that is $\partial \Sigma =\emptyset$, $S_{q^{1/2}}(\Sigma)$ is just the (Kauffman bracket) skein algebra. For the same reason as the skein module, we  denote the skein algebra as  $S_{q}(\Sigma)$.

\def \cg {\mathcal{G}}
\def \cp {\mathcal{P}}
\def \gp {(\mathcal{G},\mathcal{P})}

A {\bf marked surface} is a pair $(\mathcal{G},\mathcal{P})$, where $\cg$ is a compact surface and 
$\cp$ is a finite subset of $\partial\cg$. The boundary component that has no intersection with $\cp$ is called {\bf unmarked boundary component}.
For such a marked surface, we can define a marked three manifold 
$\MN$ where $M = \cg\times [0,1]$ and $\mathcal{N} = \cup_{x\in\cp}(\{x\}\times(0,1))$ and the orientation of $\mathcal{N}$ is given by the positive direction of $(0,1)$. Then $\MN$ is called the thickening of $\gp$, and we define $\S\gp = \S\MN$. Similarly as the pb surface, $\S\gp$ has an algebra structure given by stacking the stated tangles. When $\mathcal{P} = \emptyset$, we use $S_q(\cg)$ to denote $\S\gp$.

For an oriented surface $\Sigma$, the orientation of $\partial\Sigma$ induced from the orientation of $\Sigma$ is called the positive orientation of $\partial\Sigma$.
For any marked surface $\gp$, we can define a pb surface $\Sigma$ by the following way: (1) Let $\mathcal{Q}$ be obtained from $\cp$ be slightly moving $\cp$ along $\partial\cg$ in the positive direction of $\partial\cg$. (2) Let 
$\cg^{'}$ be obtained from $\cg$ by replacing all unmarked boundary components with interior punctures. (3) Define
$\Sigma = \cg^{'} - \mathcal{Q}$. Clearly we can also get $\gp$ from $\Sigma$.
It is easy to see  $\gp$ and $\Sigma$ have the same skein theory. 
%

Pb surface is convenient for defining triangulation. Marked surface also has it's own convenience because the thickening of the marked surface is a compact marked three manifold.


\def \CSM {S_1(M,\mathcal{N})}

\def \SMN {S_{q^{1/2}}(M,\mathcal{N})^{(N)}}

\subsection{Chebyshev polynomial and Frobenius map}\label{sub2.2}
 Chebyshev polynomials are defined by the following recurrenc relation:
$$Q_n(x) = xQ_{n-1}(x) - Q_{n-2}(x).$$
The Chebyshev polynomial of the first kind, denoted as $T_n(x)$, is defined by setting $T_0(x) = 2, T_1(x) = x$.
The Chebyshev polynomial of the second kind, denoted as $S_n(x)$, is defined by setting $S_0(x) = 1, S_1(x) = x$.

Then for $n\geq 2$, we have $T_n(x) = S_n(x) - S_{n-2}(x).$ For any two positive integers $n_1,n_2$, we have 
\be\label{Cheeq}
T_{n_1}(T_{n_2}(x)) = T_{n_1n_2}(x).
\ee

Let $(M,\mathcal{N})$ be a marked three manifold. Suppose $\beta$ is a  framed knot or stated framed arc  in $(M,\mathcal{N})$, and 
$k$ is a positive integer. We define $\beta^{(k)}$  to be obtained from $\beta$ by taking $k$ parallel copies of $\beta$ in the framing direction. 
Suppose $P(x) = \sum_{0\leq k\leq m}a_k x^k\in\mathbb{C}[x]$, then we define
$P(\alpha) = \sum_{0\leq k\leq m}a_k \alpha^{(k)} \in\SMQ$.

Suppose $q^{\frac{1}{2}}$ is a root of unity of order $N$ with $N$ being odd.
Then there is a linear map, called Frobenius map, $\cF:\CSM\rightarrow \SMQ$ \cite{bloomquist2020chebyshev,bonahon2016representations,wang2023stated}.
Let $l$ be any stated $(M,\mathcal{N})$-tangle, suppose 
$l = \alpha_1\cup\cdots\cup \alpha_k \cup\beta_1\cup\cdots\cup \beta_m$
where $\alpha_i,1\leq i\leq k,$ are stated framed arcs and $\beta_j,1\leq j\leq m,$ are framed knots.
Then 
\begin{equation}\label{Frobe}
\cF(l) = \alpha_1^{(N)}\cup\cdots\cup \alpha_k^{(N)} \cup T_N(\beta_1)\cup\cdots\cup T_N(\beta_m).
\end{equation}


Let $\alpha$ be any stated $(M,\mathcal{N})$-tangle, then $\cF(\alpha)$ is transparent \cite{bloomquist2020chebyshev,wang2023stated}, in a sense that we have the following two relations in $\SMQ$:
\begin{equation}\label{well_de}
\raisebox{-.20in}{
\begin{tikzpicture}
\filldraw[draw=white,fill=gray!20] (-0,-0.2) rectangle (1, 1.2);
\draw [line width =1pt](0.6,0.6)--(1,1);
\draw [line width =1pt](0.6,0.4)--(1,0);
\draw[line width =1pt] (0,0)--(0.4,0.4);
\draw[line width =1pt] (0,1)--(0.4,0.6);
\draw[line width =1pt] (0.6,0.6)--(0.4,0.4);
\node [right] at (1,1) {$\cF(\alpha)$};
\end{tikzpicture}
}=
\raisebox{-.20in}{
\begin{tikzpicture}
\filldraw[draw=white,fill=gray!20] (-0,-0.2) rectangle (1, 1.2);
\draw [line width =1pt](0.6,0.6)--(1,1);
\draw [line width =1pt](0.6,0.4)--(1,0);
\draw[line width =1pt] (0,0)--(0.4,0.4);
\draw[line width =1pt] (0,1)--(0.4,0.6);
\draw[line width =1pt] (0.4,0.6)--(0.6,0.4);
\node [right] at (1,1) {$\cF(\alpha)$};
\end{tikzpicture}
},
\raisebox{-.35in}{
\begin{tikzpicture}
\filldraw[draw=white,fill=gray!20] (-0.2,-0) rectangle (1.2, 1);
\draw [line width =1pt](0,1)--(0.4,0.2);
\draw [line width =1pt](1,1)--(0.6,0.2);
\draw [line width =1pt](0.5,0)--(0.6,0.2);
\filldraw[draw=black,fill=black] (0.5,-0) circle (0.09);
\node at (0.3,-0.15) {$i$};
\node at (0.7,-0.15) {\small $j$};
\node [above] at (1,1) {$\cF(\alpha)$};
\end{tikzpicture}
}=
\raisebox{-.35in}{
\begin{tikzpicture}
\filldraw[draw=white,fill=gray!20] (-0.2,-0) rectangle (1.2, 1);
\draw [line width =1pt](0,1)--(0.4,0.2);
\draw [line width =1pt](1,1)--(0.6,0.2);
\draw [line width =1pt](0.5,0)--(0.4,0.2);
\filldraw[draw=black,fill=black] (0.5,-0) circle (0.09);
\node at (0.3,-0.15) {$i$};
\node at (0.7,-0.15) {\small $j$};
\node [above] at (1,1) {$\cF(\alpha)$};
\end{tikzpicture}
} 
\end{equation}
where $i,j= -,+$.

\def \SMO {S_1(M,N)}

Then we try to define an action of $\SMO$ on $\SMQ$.
For any elements $\alpha = c_1\alpha_1+\cdots + c_m\alpha_m\in\SMO$ and $l =k_1 l_1+\cdots + k_nl_n\in \SMQ$ where $c_j,k_i\in\mathbb{C},$ $\alpha_j$ and $l_i$ are stated tangles in $\MN$ and $\alpha_j\cap l_i=\emptyset$,
  we define 
\begin{equation}\label{define}
\alpha\cdot l=
\cF(\alpha)\cup l =\sum_{1\leq j\leq m,1\leq i\leq n} c_j k_i \cF(\alpha_j)\cup l_i.
\end{equation}
From equation \eqref{well_de} and the fact that $\cF$ preserves  relations \eqref{cross}-\eqref{hight}, the above action is well-defined.



\def \hF {\hat{\cF}}

We use $\SMN$ to denote Im$\cF$. 
\bt\label{multi}
$\SMN$ has a commutative algebra structure, which makes $\cF:\SMO\rightarrow \SMN$ a surjective algebra homomorphism.
\et
\begin{proof}
For any two elements $\alpha,l\in \SMO$, first we prove $\alpha\cdot\cF( l) = l\cdot\cF(\alpha)$. Suppose 
$\alpha = c_1\alpha_1+\cdots c_m\alpha_m\in\SMO$ and $l =k_1 l_1+\cdots + k_nl_n\in \SMO$ where $c_j,k_i\in\mathbb{C},$ $\alpha_j$ and $l_i$ are stated tangles in $\MN$ and $\alpha_j\cap l=\emptyset$.
Then from equation \eqref{define}, we know 
\begin{align}\label{surj}
\alpha\cdot\cF( l)=\sum_{1\leq j\leq m,1\leq i\leq n} c_j k_i \alpha_j\cdot \cF(l_i)
 =\sum_{1\leq j\leq m,1\leq i\leq n} c_j k_i \cF(\alpha_j)\cup\cF( l_i) = l\cdot\cF(\alpha).
\end{align}

For any two elements $x,y\in \SMN$, suppose $x= \cF(u)$ where $u\in\SMO$. Then we define $xy = u\cdot y$. To show this multiplication is well-defined, we have to show $xy$ is independent of the choice of $u$. Suppose we also have $x= \cF(u^{'})$ where $u^{'}\in\SMO$. We also suppose $y = \cF(v)$ where $v\in\SMO$. Then 
$$ u\cdot y =u\cdot \cF(v) = v\cdot \cF(u) = v\cdot \cF(u^{'}) = u^{'}\cdot\cF(v) = u^{'}\cdot y.$$
We also have 
$$xy =  u\cdot y =u\cdot \cF(v) = v\cdot \cF(u)=v\cdot x = yx.$$ Thus the above defined algebra structure is commutative.

Equation \eqref{surj} indicates $\cF:\SMO\rightarrow \SMN$ is an algebra homomorphism.
\end{proof}
%
%


\def\SqS {S_{q^{1/2}}(\Sigma)}
\def\CSq {S_{1}(\Sigma)}
\def\SQN {S_{q^{1/2}}(\Sigma)^{(N)}}

For a pb surface $\Sigma$, we have $\cF : \CSq\rightarrow \SqS$ is an algebra homomorphism. We also use 
$\SQN$ to denote Im$\cF$. Then the multiplication structure of $\SQN$ inherited from $\SqS$ coincides with the one defined in Theorem \ref{multi}.
 Then we have $\cF : \CSq\rightarrow \SQN$ is an isomorphism and $\SQN\subset C(\SqS)$
\cite{bloomquist2020chebyshev,korinman2019classical,wang2023stated}.
Obviously $\SqS$ has an $\SQN$-module structure defined by multiplication. We know $\SqS$ also has 
an $\CSq$-module structure. Then these two module structures are equivalent via $Id_{\SqS}$ and $\cF$.
That is, for any  $x\in\CSq,y\in\SqS$, we have $x\cdot y= \cF(x)y$.

%


%

\def \fs{f_{\sharp}}

\def \Sqo{S_q(M_1)}
\def \Sqt{S_q(M_2)}
\def \SqoN{S_q(M_1)^{(N)}}
\def \SqtN{S_q(M_2)^{(N)}}

\section{On skein modules}\label{section3}
The  skein module $\Sq$ for a three manifold $M$ does not involve boundary and is only defined by relations \eqref{cross} and \eqref{unknot}.
%
%
Suppose $q^2$ is a primitive $N$-th root of unity with $N$ being odd. Then there exists a linear map
$$\cF: S_{\zeta}(M)\rightarrow S_q(M),$$
where $\zeta = q^{N} = \pm1$ \cite{bloomquist2020chebyshev,bonahon2016representations}.
The definition here is the same with equation \eqref{Frobe}, except $l$ does not contain arcs.

We still use $\SqN$ to denote Im$\cF$.
From \cite{bloomquist2020chebyshev}, we know elements in $\SqN$ are transparent. Similarly as in subsection \ref{sub2.2}, $\Sq$ has an $S_{\zeta}(M)$-module structure, and
$\SqN$ has a commutative algebra structure, which makes $\cF:\Sz\rightarrow \SqN$ a surjective algebra homomorphism. For any  $\eta\in\Sz$ and  $\gamma\in \Sq$,  we use $\eta\cdot \gamma$ to denote the action of $\Sz$ on $\Sq$ (the same notation as in subsection \ref{sub2.2}).  

For a closed pb surface $\Sigma$, we  have $\cF:S_{\zeta}(\Sigma)\rightarrow S_q(\Sigma)$ is an algebraic injection \cite{bloomquist2020chebyshev,bonahon2016representations}.
We use $ S_q(\Sigma)^{(N)}$ to denote $\text{Im}\cF$, then $ S_q(\Sigma)^{(N)}$ also acts on $S_q(\Sigma)$ by the multiplication. And for all $x\in S_{\zeta}(\Sigma),y\in S_q(\Sigma)$, we have  $x\cdot y = \cF(x)y$.



In this section, we will always assume $q^2$ is a primitive $N$-th root of unity with $N$ being odd. 
\begin{theorem}(\cite{talk})\label{skein}
For any compact  three manifold $M$, we have $\Sq$ is finitely generated over $\Sz$.
\end{theorem}

\def \SqN {S_{\zeta}(M)}

In this section,
we use the connection between Theorem \ref{skein} and the finiteness of the skein module of the three manifold $M$ over $S_q(\partial M)$ to give sufficient conditions for $\Sq$ being finitely generated over $S_q(\partial M)$. We also provide a way to estimate an upper bound for $\dim_{\SqN}\Sq$ when $M$ is compact, which will be used in subsection \ref{sub6.1} to give a precise upper bound for $\dim_{\SqN}\Sq$.


%
%


\subsection{Some functorialities}

Let $f: M_1\rightarrow M_2$ be an embedding  between two three manifolds. We know $f$ induces a linear map 
$\fs : S_q(M_1)\rightarrow S_q(M_2)$ and an algebra homomorphism $\fs:S_{\zeta}(M_1)\rightarrow S_{\zeta}(M_2)$. Actually we have the following lemma.

\bl\label{functor}
Let $f: M_1\rightarrow M_2$ be an embedding between two three manifolds.
Then $\fs$ restricts to an algebra homomorphism $\fs:\SqoN\rightarrow \SqtN,$ and we have the following commutative diagram:
$$
\begin{tikzcd}
S_{\zeta}(M_1)  \arrow[r, "\fs"]
\arrow[d, "\cF"]  
&  S_{\zeta}(M_2)  \arrow[d, "\cF"] \\
 S_q(M_1)  \arrow[r, "\fs"] 
&  S_q(M_2)\\
\end{tikzcd}.
$$
What's more, $\fs$ respects the module structures in a sense that, for any $$x\in S_{\zeta}(M_1),y\in\Sqo,$$ we have
$\fs(x\cdot y) = \fs(x)\cdot \fs(y).$
\el

\bl\label{keylemma}
Let $f: M_1\rightarrow M_2$ be an embedding between two three manifolds. Suppose 
$\fs : S_q(M_1)\rightarrow S_q(M_2)$ is surjective and $S_q(M_1)$ is finitely generated over $S_{\zeta}(M_1)$. Then 
we have $S_q(M_2)$ is finitely generated over $S_{\zeta}(M_2)$ and $$\text{dim}_{S_{\zeta}(M_2)}\Sqt\leq \text{dim}_{S_{\zeta}(M_1)}\Sqo.$$
\el
\bpr
Suppose $\text{dim}_{S_{\zeta}(M_1)}\Sqo = n$ and $\Sqo$ is  generated by $x_1,\cdots,x_n$ over $S_{\zeta}(M_1)$. Then 
$\Sqt$ is finitely generated by $\fs(x_1),\cdots,\fs(x_n)$ over $S_{\zeta}(M_2)$ because of the surjectivity of $\fs$  and
Lemma \ref{functor}. 
\epr

%

\subsection{Boundary action on skein modules}

Let $M$ be a three manifold. Then the skein algebra $\Sqp$ has an action on $\Sq$. We can identity $\partial M\times [0,1]$ with
a regular closed tubular neighborhood $U(\pM)$ of $\partial M$ such that $\pM$ is $\pM\times \{1\}$. We use 
$L$ to denote the embedding from $U(\pM)$ to $M$. 
Then for any skein $l\in\Sqp$, we have $\Ls(l)$ is the skein $l$ in $\Sq$.
Then for any skein $\alpha$ in $\Sqp$ and $\beta\in \Sq$, we first isotope $\beta$ such that $\beta\in
M-U(\pM)$, then define $\alpha*\beta \define \Ls(\alpha)\cup \beta$.


We know that $\SqN$ also acts on $\Sq$. Then
it is easy to show these two actions commute with each other, that is, for any elements
$\alpha\in\Sqp, \beta\in \SqN,\eta\in\Sq$ we have $\alpha*(\beta\cdot\eta) = \beta\cdot (\alpha*\eta)$.

For elements $l_1\in S_{\zeta}(\partial M), l_2\in \Sqp, l_3\in \Sq$, we  have 
\begin{align}\label{666}
(l_1\cdot l_2)*l_3 = (\cF(l_1)l_2)*l_3 =\cF( l_1)*(l_2*l_3) = \Ls(l_1)\cdot(l_2*l_3).
\end{align}

\subsection{Relations between $\Sq$ over $\SqN$ and $\Sq$ over $S_q(\partial M)$}

When the ground ring is $\mathbb{Q}(q)$, Detcherry conjectured that the skein module of a compact three manifold $M$ is finitely generated as a module over the skein algebra  $S_q(\partial M)$  \cite{detcherry2021infinite}. When the ground ring is $\mathbb{C}$ and $q^2$ is a root of unity of odd order, this conjecture fails. A simple counterexample is $S^1\times S^2$. But Detcherry's conjecture still holds for some   three manifolds under this section's context , that is when the ground ring is $\mathbb{C}$ and $q^2$ is a root of unity of odd order, for example the lens spaces and the complement of two bridge knot (or link) \cite{bullock2005kauffman,hoste19932,le2014kauffman, le2006colored}. In this subsection, we will give sufficient conditions for $\Sq$ being finitely generated over $S_q(\partial M)$.

\def \dim{\text{dim}}

\bt\label{finitesurface} (Theorem 3.11 in \cite{abdiel2017localized})
Suppose $\Sigma$ is a finite type surface. Then $\Si$ is finitely generated over $S_{\zeta}(\Sigma)$ (or $\SiN$).
\et

We state a convention that $S_q(\emptyset)$ is $\mathbb{C}$.

\bp\label{relation1}
  Let $M$ be  a three manifold with $\partial M$ being a finite type surface.
 Suppose $\Sq$ is finitely generated over $\Sqp$. Then $\Sq$ is finitely generated over $\SqN$, and we have $$\text{dim}_{\SqN}\Sq\leq (\text{dim}_{\Sqp}\Sq)
(\text{dim}_{S_{\zeta}(\partial M)}\Sqp). $$
\ep
\bpr
When $\partial M=\emptyset$, obviously we have 
$$\text{dim}_{\SqN}\Sq\leq \text{dim}_{\mathbb{C}}\Sq. $$

Then we look at the case when $\partial M\neq \emptyset.$
Suppose $\text{dim}_{\Sqp}\Sq = r$, and $\Sq$ is generated by $A_1,\dots, A_r$ as an $\Sqp$-module. 
From Theorem \ref{finitesurface}, we can
suppose $\text{dim}_{S_{\zeta}(\partial M)}\Sqp= s$, and $\Sqp$ is generated by $B_1,\dots, B_s$ as an $S_{\zeta}(\partial M)$-module. 
It suffices to show $\Sq$ is generated by $B_{i}* A_{j},1\leq i\leq s, 1\leq j\leq r$ as an $\SqN$-module.

Let $C$ be any element in $\Sq$. Then $C = C_1 * A_1+\dots+ C_r * A_r$ where $C_1,\dots, C_r\in \Sqp$. For each $1\leq j\leq r$, we can suppose
$$C_j = \sum_{1\leq i\leq s} D_{j,i}\cdot B_i$$
where $D_{j,i}\in S_{\zeta}(\partial M)$.
Then we have 
\begin{align*}
C &= \sum_{1\leq j\leq r}C_j*A_j = \sum_{1\leq j\leq r}(\sum_{1\leq i\leq s} D_{j,i}\cdot B_i)*A_j\\&=
\sum_{1\leq j\leq r,1\leq i\leq s}(D_{j,i}\cdot B_i)*A_j =  \sum_{1\leq j\leq r,1\leq i\leq s}\Ls(D_{j,i})\cdot(B_i*A_j).
\end{align*}
where the last equality is because of equation \eqref{666}, and $\Ls(D_{j,i})\in\SqN$, for $1\leq j\leq r,1\leq i\leq s$.
\epr

%

\bp\label{relation2}
Let $M$ be a three manifold. If $\Sq$ is finitely generated over $\SqN$ and $\Sz$ is finitely generated over $\Szp$, then $\Sq$ is finitely generated over $\Sqp$ and $\text{dim}_{\Sqp}\Sq\leq (\text{dim}_{\Szp}\Sz)
(\text{dim}_{\SqN}\Sq$).
\ep
\bpr

Assume $\text{dim}_{\SqN}\Sq = n$, and $\Sq$ is generated by $b_1,\dots, b_n$ over $\SqN$. 
Assume $\text{dim}_{\Szp}\Sz= m$, and $\Sz$ is generated by $c_1,\dots, c_m$ over $\Szp$. 
It suffices to show $\Sq$ is generated by $c_{i}\cdot b_{j},1\leq i\leq m, 1\leq j\leq n$ over $\Sqp$.

For any $a\in\Sq$,
we can suppose $a = \sum_{1\leq i\leq n} a_i\cdot b_i$ where $a_i\in\Sz,1\leq i\leq n$.
For each $1\leq i\leq n$, we suppose $a_i = \sum_{1\leq j\leq m} d_{i,j}* c_j$ where $d_{i,j}\in \Szp,1\leq i\leq n,1\leq j\leq m$. Then we have
\begin{align*}
a& = \sum_{1\leq i\leq n} a_i\cdot b_i  = \sum_{1\leq i\leq n,1\leq j\leq m} (d_{i,j}* c_j)\cdot b_i
=\sum_{1\leq i\leq n,1\leq j\leq m}  \cF(d_{i,j})*(c_j\cdot b_i) 
\end{align*}
where $\cF(d_{i,j})\in \Sqp,1\leq i\leq n,1\leq j\leq m$.
\epr

\bt\label{thm3.9}
Let $M$ be a compact   three manifold. If $\Sz$ is finitely generated over $\Szp$, then $\Sq$ is finitely generated over $\Sqp$ and $$\text{dim}_{\Sqp}\Sq\leq (\text{dim}_{\Szp}\Sz)
(\text{dim}_{\SqN}\Sq).$$
\et
\bpr
Theorem \ref{skein} and Proposition \ref{relation2}.
\epr

%

\subsection{Handlebody and compression body}
A three manifold $H$ is called a {\bf compression body} if there exists a connected totally
closed surface $\Sigma$ such that $H$ is obtained from $\Sigma\times[0,1]$ by attaching 2-handles
along mutually disjoint loops in $\Sigma\times \{1\}$ and filling in some resulting 2-sphere
boundary components with three dimensional solid balls. We denote $\Sigma \times \{0\}$ by $\partial^{+} H$ 
and $\partial H - \partial^{+} H$ by $\partial ^{-}H$. We define the genus of $H$ to be the genus of $\Sigma$.
 A compression body $H$ is called a {\bf handlebody} if $\partial^{-} H =\emptyset$ \cite{saito2005lecture}.

\bt(\cite{saito2005lecture})\label{splitting}
Let $M$ be a connected compact  three manifold. Then $M = H\cup H^{'}$ where $H$ is a handlebody and 
$H^{'}$ is a compression body such that $\partial^{-}H^{'} = \partial M$ and $\partial^{+}H^{'} = \partial H$.

\et

The decomposition in Theorem \ref{splitting} is called the Heegaard splitting of $M$. We define the genus (or Heegaard genus) of $M$
 to be the minimum  genus of the handlebody among all the Heegaard splittings.

Let $(M,\mathcal{N})$ be a compact marked three manifold,  we define the genus of $(M,\mathcal{N})$ to be the genus of $M$. 
When $\mathcal{N}$ is empty we could just call the genus of $(M,\emptyset)$ as the genus of $M$.

For any nonnegative integer $k$, we $P_k$ to denote the surface obtained from $S^{2}$ by removing $k$ points.
Let $H$ be a handlebody of genus $g$, then $H$ has the same skein theory with $P_{g+1}\times [0,1]$. 
\subsection{An upper bound for $\dim_{\SqN}\Sq$}
Let $M$ be a compact   three manifold. We can assume $M$ is a connected compact  three manifold with genus $g$.  Then
suppose the Heegaard splitting of $M$ is $M= H\cup H^{'}$ where $H$ is a handlebody of genus $g$ and $H^{'}$ is a compression body.  Then we can isotope
any skein $\alpha$ in $M$  into $H$. Thus the embedding $f: H\rightarrow M$ induces a surjective linear map $\fs:S_q(H)\rightarrow\Sq$. From Lemma \ref{keylemma},  we can get 
\begin{align}\label{eq_handle}
\dim_{\SqN}\Sq\leq \dim_{S_{\zeta}(H)}S_q(H) = \dim_{S_{\zeta}( P_{g+1})}S_q( P_{g+1})= \dim_{S_{q}( P_{g+1})^{(N)}}S_q( P_{g+1}).
\end{align}


From Theorem \ref{finitesurface}, we know
$\dim_{S_q( P_{g+1})^{(N)}}S_q( P_{g+1})$ is finite. Thus
the above discussion also offers a proof for Theorem \ref{skein}.

\def \xM {\mathcal{X}(M)}

\subsection{Finiteness of the reduced skein module}\label{sub3.4}
Let $M$ be a three manifold. We use $\mathcal{X}(M)$ to denote $ Hom(\pi_1(M),SL_2(\mathbb{C}))/\simeq$
where $Hom(\pi_1(M),SL_2(\mathbb{C}))$ is the set of group homomorphisms from $\pi_1(M)$ to 
$SL_2(\mathbb{C})$, and  two homomorphisms $\rho_1\simeq\rho_2$ if and only if $\text{Trace}(\rho_1(x))
=\text{Trace}(\rho_2(x))$ for all $x\in\pi_1(M)$. Then for any $\rho\in\mathcal{X}(M)$, it induces an algebra homomorphism $f_{\rho}:\Sz\rightarrow \mathbb{C}$, furtermore this correspondence is a bijection from 
$\mathcal{X}(M)$ to the set of algebra homomorphisms from $\Sz$ to $\mathbb{C}$
\cite{przytycki1997skein}.

Let $I$ be an ideal of $\SqN$,
 we  
$$I\Sq = \{x_1 \cdot y_1+\cdots+ x_k\cdot y_k\mid k\in\mathbb{N}, x_i\in \SqN, y_i\in\Sq, 1\leq i\leq k\}$$  
(note that $x_1 \cdot y_1+\cdots+ x_k\cdot y_k $ is defined to be $0$ when $k=0$).
Then $I\Sq$ is an $\SqN$-submodule of $\Sq$.

 For any element $\rho\in\xM$, the {\bf reduced skein module} $S_q(M)_{\rho}$ of $M$ with respect to $\rho$ is defined to be 
\begin{equation}\label{def_reduced}
S_q(M)/I_{\rho}
\end{equation}
where $I_{\rho} = (\text{Ker}f_{\rho})\Sq$ 
\cite{detcherry2023kauffman,talk}. 



\begin{remark}\label{rem3.6}
For any maximal ideal $I$ of $\SqN$, which is point in $\mathcal{X}(M)$,
the reduced skein module is still an $\SqN$-module since $I\Sq$ is an $\SqN$-submodule. Actually 
$\SqN$-module structure for $\Sq_{I}$ coincides with the linear structure over $\mathbb{C}$ by the identity map.
The maximal ideal $I$ of $\SqN$ corresponds to an algebra homomorphism $f_I: \SqN\rightarrow \mathbb{C}$, then
for any element $\alpha\in\SqN,\beta\in\Sq$, we have 
$$\alpha\cdot(\beta+I\Sq) = \alpha\cdot\beta+ I\Sq = f_I(\alpha)(\beta+I\Sq).$$
\end{remark}

From Remark \ref{rem3.6}, we know
 the reduced skein module is obtained by sending elements in $\SqN$ to  complex numbers, we can easily get the following result by using Theorem \ref{skein}.
\bt\label{reduced_dim}
Let $M$ be a compact  three manifold, and $\rho\in \xM$. Then 
$S_q(M)_{\rho}$ is finite dimensional, and $\dim_{\mathbb{C}}\Sq_{\rho}\leq \dim_{\SqN}\Sq$.
\et

%
%
%


\def \SMN {S_1(M,\mathcal{N})}

\section{On stated skein modules}\label{section4}
In this section $q^{\frac{1}{2}}$ is always a root unicty of order $N$ with $N$ being odd. 
The main goal of this section is to prove Theorem \ref{1.1}.
Here we restate the Theorem for reader's  convenience.  
\begin{theorem}\label{th1.1}
 $\SMQ$ is finitely generated over $\SMN$ for any compact marked three manifold $\MN$. 
\end{theorem}

\subsection{Some functorialities}

Let $f: (M_1,\mathcal{N}_1)\rightarrow (M_2,\mathcal{N}_2)$ be an embedding  between two marked three manifolds. Then $f$ induces a linear map 
$\fs : S_{q^{1/2}}(M_1,\mathcal{N}_1)\rightarrow S_{q^{1/2}}(M_2,\mathcal{N}_2)$
, which restricts to an algebra homomorphism $\fs:S_{q^{1/2}}(M_1,\mathcal{N}_1)^{(N)}\rightarrow S_{q^{1/2}}(M_2,\mathcal{N}_2)^{(N)}$. Similarly as in Lemma \ref{functor}, we have 
  $f$ also induces an algebra homomorphism $\fs :S_1(M_1,\mathcal{N}_1)\rightarrow S_1(M_2,\mathcal{N}_2)$,
 and $\fs$ respects the module structures, in a sense that, for any $x\in S_1(M_1,\mathcal{N}_1), y\in S_{q^{1/2}}(M_1,\mathcal{N}_{1})$ we have $\fs(x\cdot y) = \fs(x)\cdot \fs(y)$.

\bl\label{keylemma1}
Let $f: (M_1,\mathcal{N}_1)\rightarrow (M_2,\mathcal{N}_2)$ be an embedding  between two marked three manifolds. Suppose 
$\fs :  S_{q^{1/2}}(M_1,\mathcal{N}_1)\rightarrow S_{q^{1/2}}(M_2,\mathcal{N}_2)$ is surjective and Theorem  \ref{th1.1} holds for $(M_1, \mathcal{N}_1)$. Then 
Theorem  \ref{th1.1} also holds for $(M_2, \mathcal{N}_2)$ and $$\text{dim}_{S_1(M_2,\mathcal{N}_2)}S_{q^{1/2}}(M_2,\mathcal{N}_2)\leq \text{dim}_{S_1(M_1,\mathcal{N}_1)}S_{q^{1/2}}(M_1,\mathcal{N}_1).$$
\el
\bpr
The proof here is the same with Lemma \ref{keylemma}.
\epr

\subsection{Good strict subsurface}
Let $\MN$ be a compact marked three manifold, and $\Sigma$ be a properly  embedded surface in $M$ such that 
every component of $\Sigma$ intersects $\mathcal{N}$. Let $V(\Sigma)$ be a closed regular neighborhood of $\Sigma$, that is, $V(\Sigma)$ is isomorphic to $\Sigma\times[0,1]$ by an orientation preserving diffeomorphism such that 
$\Sigma$ is identified with $\Sigma\times\{1/2\}$ and $\partial\Sigma\times[0,1]\subset \partial M$. Then define
$Cut_{\Sigma}\MN = (M^{'},\mathcal{N}^{'})$ where $M^{'} = M - \Sigma\times(0,1)$ and 
$\mathcal{N}^{'} = \mathcal{N} - \Sigma\times[0,1]$ (here we identify $\Sigma\times[0,1]$ with $V(\Sigma)$).
$\Sigma$ is called a 
{\bf good strict subsurface} of $(M,\mathcal{N})$,
if  $Cut_{\Sigma}\MN= (M^{'},\mathcal{N}^{'})$ is isomorphic to the thickening of a marked surface.

\def\SqS {S_{q^{1/2}}(\Sigma)}
\def\SQN {S_{q^{1/2}}(\Sigma)^{(N)}}

\def \gq {\S\gp}
\def \gqN {\S\gp^{(N)}}

\bt\label{statedsurface}
Let $\gp$ be any marked surface, then $\gq$ is finitely generated over $S_1(\mathcal{G},\mathcal{P})$ (or $\gqN$).
\et
\bpr
We can assume $\cg$ is connected. When $\cp$ is empty, the stated skein algebra is just the skein algebra. This case was proved in Theorem \ref{skein}. When $\cp$ is not empty, Theorem \ref{statedsurface} was proved in \cite{korinman2021unicity, wang2023stated}.
\epr

\begin{proposition}(Proposition 8.14 \cite{wang2023stated})\label{pp8.14}
Assume $(M,\mathcal{N})$ is a marked three manifold which contains a good strict subsurface $\Sigma$. Then we have $\SMQ$ is a finitely generated  $\SMN$-module, and 
$$\text{dim}_{\SMN}S_{q^{1/2}}(M,\mathcal{N})\leq \text{dim}_{S_{q^{1/2}}(Sl_{\Sigma}(M,\mathcal{N}))^{(N)}}S_{q^{1/2}}(Sl_{\Sigma}(M,\mathcal{N})).$$

\end{proposition}
\bpr
Here we briefly state the proof in \cite{wang2023stated}. We use $(M^{'},\mathcal{N}^{'})$ to denote $Sl_{\Sigma}(M,\mathcal{N})$. Since $\Sigma$ is a good strict subsurface, then $(M^{'},\mathcal{N}^{'})$ has the same skein theory with the thickening of a  marked surface. Then Theorem \ref{statedsurface} indicates $S_{q^{1/2}}(M^{'},\mathcal{N}^{'})$ is finitely generated over $S_{1}(M^{'},\mathcal{N}^{'})$.

Meanwhile  the obvious embedding $f:(M^{'}, N^{'})\rightarrow (M,N)$ induces a surjective linear map
$f_{\sharp}:S_{q^{1/2}}(M^{'}, N^{'})\rightarrow S_{q^{1/2}}(M,N)$. Then Lemma \ref{keylemma1} completes the proof. And we can get 
$$\text{dim}_{S_{1}(M,\mathcal{N})}S_{q^{1/2}}(M,\mathcal{N})\leq \text{dim}_{S_{1}(M^{'},\mathcal{N}^{'})}S_{q^{1/2}}(M^{'},\mathcal{N}^{'}).$$
Since $(M^{'},\mathcal{N}^{'})$ is isomorphic to the thickening of a marked surface,  we have 
$$\text{dim}_{S_{1}(M^{'},\mathcal{N}^{'})}S_{q^{1/2}}(M^{'},\mathcal{N}^{'})=\text{dim}_{S_{q^{1/2}}(Sl_{\Sigma}(M,\mathcal{N}))^{(N)}}S_{q^{1/2}}(Sl_{\Sigma}(M,\mathcal{N})).$$
\epr

For any two nonnegative integers $g,p$, we use $\Sigma_{g,p}$ to denote the compact oriented surface with $p$ boundary components and genus $g$. 
For any two nonnegative integers $a,b$ with $b>0$, we use ${}_{a,b}\Sigma$ to denote the pb surface obtained from $\Sigma_{a,1}$ by removing $b$ points on the unique boundary component of $\Sigma_{a,1}$, we use ${}_{a,b}\cg$ to denote the marked surface $(\Sigma_{a,1},\cp)$ where $\cp$ is the subset of  the unique boundary component of $\Sigma_{a,1}$ consisting of $k$ points.
For the positive integer $k$, we use $D_k$ to denote the pb surface obtained from the disk by removing $k$ points on the boundary, we use $\cg_k$ to denote the marked surface $(D,K)$ where $D$ is the disk and $K\subset\partial D$ consisting of $k$ points. We call $D_2$ the  bigon, and  use $\lambda_{q^{1/2}}$ to denote 
$\text{dim}_{S_{q^{1/2}}(D_2)^{(N)}}S_{q^{1/2}}(D_2)$.

Let $k,g$ be two nonnegative integers, we use $M_{g,k}$ to denote the marked three manifold $(M,\mathcal{N})$, where
$M=\Sigma_{g,0}\times[0,1]$ and $N\subset \Sigma_{g,0}\times
\{1\}$ such that $\mathcal{N}$ has $k$ components. Note that $M_{g,k}$ is defined up to isomorphism.

\bl\label{compression}
For any two nonnegative integers $g,k$, we have $S_{q^{1/2}}(M_{g,k})$ is finitely generated over 
$S_{1}(M_{g,k})$. When $k>0$, we have 
$$\text{dim}_{S_{1}(M_{g,k})}S_{q^{1/2}}(M_{g,k})\leq  \text{dim}_{S_{q^{1/2}}({}_{g,k}\Sigma)^{(N)}}S_{q^{1/2}}({}_{g,k}\Sigma).$$
\el
\bpr
When $k = 0$, the marking set of $M_{g,k}$ is empty, and the stated skein module is just the skein module. 
Then Theorem \ref{skein} indicates this case.

Then we look at the case when $k$ is positive.
From the definition, we know $M_{g,k}$ is just $\Sigma_{g,0}\times [0,1]$ with $k$ markings on $\Sigma_{g,0}\times \{1\}$. We use $\mathcal{N}$ to denote the union of all the markings of $M_{g,k}$.
Let $c$ be a closed curve in $\Sigma_{g,0}$  such that $c$ bounds an embedded disk in $\Sigma_{g,0}$ and $c\times\{1\}$ intersects $\mathcal{N}$ transversely in exactly one point.
Then $c\times [0,1]$ is a strict subsurface of $M_{g,k}$. We also have 
$Sl_{c\times [0,1]}(M_{g,k})$  is isomorphic with the thickening of ${}_{g,k}\cg\cup \cg_{1}$.
Thus $c\times [0,1]$ is a good strict subsurface of $M_{g,k}$. Then Theorem \ref{statedsurface} and
Proposition \ref{pp8.14} indicates 
$S_{q^{1/2}}(M_{g,k})$ is finitely generated over 
$S_{1}(M_{g,k})$, and 
$$\text{dim}_{S_{1}(M_{g,k})}S_{q^{1/2}}(M_{g,1})\leq\text{dim}_{S_{q^{1/2}}({}_{g,k}\cg)^{(N)}}S_{q^{1/2}}({}_{g,k}\cg)= \text{dim}_{S_{q^{1/2}}({}_{g,k}\Sigma)^{(N)}}S_{q^{1/2}}({}_{g,k}\Sigma).$$

\epr

\subsection{Attaching 2-handles and cutting out  open 3-balls for marked three manifolds}
Let $(M,\mathcal{N})$ be a marked three manifold, and $c_1,\cdots,c_k$ be a collection of closed disjoint curves on $\partial M$ such that for each $1\leq i\leq k$ there is no intersection between $c_i$ and $\mathcal{N}$.
Then we can define a new marked three manifold $(M^{'}, \mathcal{N}^{'})$ where $M^{'}$ is obtained from $M$
by attaching 2-handles along $c_1,\cdots,c_k$ and $\mathcal{N}^{'} = \mathcal{N}$. We say $(M^{'},\mathcal{N}^{'})$ is obtained from $(M,\mathcal{N})$ by attaching 2-handles along $c_1,\cdots,c_k$. We use $f:(M,\mathcal{N})\rightarrow (M^{'},\mathcal{N}^{'})$ to denote the obviously embedding.
 Then we have the following obvious Lemma.

\bl\label{2_handle}
$f:(M,\mathcal{N})\rightarrow (M^{'},\mathcal{N}^{'})$ induces a surjective linear map 
$$\fs:S_{q^{1/2}}(M,\mathcal{N})\rightarrow S_{q^{1/2}}(M^{'},\mathcal{N}^{'}).$$
\el

\bl\label{handle}
Let $(M,\mathcal{N})$ be a marked three manifold such that $\SMQ$ is finitely generated over $\SMN$. Suppose $(M^{'},\mathcal{N}^{'})$
is obtained from $(M,\mathcal{N})$ by attaching 2-handles.
Then $S_{q^{1/2}}(M^{'},\mathcal{N}^{'})$ is also finitely generated over 
$S_{1}(M^{'},\mathcal{N}^{'})$, and 
$$\text{dim}_{S_{1}(M^{'},\mathcal{N}^{'})}S_{q^{1/2}}(M^{'},\mathcal{N}^{'})\leq \text{dim}_{S_{q^{1/2}}(M,\mathcal{N})^{(N)}}S_{q^{1/2}}(M,\mathcal{N}).$$
\el
\bpr
Lemmas \ref{keylemma1} and \ref{2_handle}.
\epr

\begin{remark}\label{ball}
Let $(M,\mathcal{N})$ be a marked three manifold, and $B$ be an embedded three dimensional solid ball contained in the interior of $M$. Let $M^{'}$ be obtained from $M$ by cutting out the interior of $B$, conversely $M$ is obtained from $M^{'}$ by filling $\partial B$ with a solid three dimensional ball. Clearly, $(M,\mathcal{N})$ and $(M^{'},\mathcal{N})$ have the same skein theory.
\end{remark}

\bl\label{keyl}
Let $(M,\mathcal{N})$ be a marked three manifold such that $M$ is a compression body $H$ and $\mathcal{N}$ is a nonempty subset of $ \partial^{-}H$. Then 
$\SMQ$ is finitely generated over $\SMN$, and 
$$\text{dim}_{\SMN}S_{q^{1/2}}(M,\mathcal{N})\leq \text{dim}_{S_{q^{1/2}}({}_{g,k}\Sigma)^{(N)}}S_{q^{1/2}}({}_{g,k}\Sigma)$$
where $g$ is the genus of $H$ and $k$ is the number of components of $\mathcal{N}$.
\el
\bpr
From the definition of the compression body, we know $(M,\mathcal{N})$ is obtained from $M_{g,k}$ by attaching 2-handles
along mutually disjoint loops in $\Sigma\times \{1\}$ and filling in some resulting 2-sphere
boundary components with three dimensional solid balls. Then Lemmas \ref{compression}, \ref{handle} and Remark \ref{ball}
complete the proof. 
\epr

\subsection{Proof for Theorem \ref{th1.1}}\label{sub666}
 Let $(M,\mathcal{N})$ be a compact marked three manifold.  We can  assume $(M,\mathcal{N})$ is a connected compact marked three manifold. When $\mathcal{N} = \emptyset$, Theorem \ref{skein} indicates this case. Then suppose $\mathcal{N}\neq\emptyset$. From Theorem \ref{splitting}, we have  $M = H\cup H^{'}$ where $H$ is a handlebody and 
$H^{'}$ is a compression body such that $\partial^{-}H^{'} = \partial M$ and $\partial^{+}H^{'} = \partial H$.
Then $\mathcal{N}\subset \partial M = \partial^{-}H^{'}$. Let $f$ be the obvious embedding from $(H^{'}, \mathcal{N})$ to $(M,\mathcal{N})$. Then $f$ induces a surjective linear map $\fs : S_{q^{1/2}}(H^{'}, \mathcal{N}) \rightarrow S_{q^{1/2}}(M, \mathcal{N})$
since we can isotope all stated $(M,\mathcal{N})$-tangles away from $H$. Then Lemmas \ref{keylemma1} and \ref{keyl}
indicate $\SMQ$ is finitely generated over $\SMN$, and 
\begin{equation}\label{eqkey}
\text{dim}_{\SMN}S_{q^{1/2}}(M,\mathcal{N})\leq \text{dim}_{S_{q^{1/2}}({}_{g^{'},k}\Sigma)^{(N)}}S_{q^{1/2}}({}_{g^{'},k}\Sigma)
\end{equation}
where $g^{'}$ is the genus of $H^{'}$ and $k$ is the number of components of $\mathcal{N}$.

\begin{remark}\label{rem_upper}
Let $(M,\mathcal{N})$ be a connected compact marked three manifold with $\mathcal{N}\neq \emptyset$. Suppose the genus of $(M,\mathcal{N})$ is $g$ and $\mathcal{N}$ has $k$ components. Then equation \eqref{eqkey} indicates 
$$\text{dim}_{\SMN}S_{q^{1/2}}(M,\mathcal{N})\leq \text{dim}_{S_{q^{1/2}}({}_{g,k}\Sigma)^{(N)}}S_{q^{1/2}}({}_{g,k}\Sigma).$$

\end{remark}

\bt
Conjecture 8.16 in \cite{wang2023stated} holds when $n=2$.
\et

\def \XMN {\mathcal{X}\MN}

\subsection{Reduced stated skein modules}\label{sub_reduced}
Let $(M,\mathcal{N})$ be a marked three manifold. 
For any ideal $I$ of $\SMN$, we define 
$$I\SMQ= \{x_1\cdot\alpha_1+\cdots+x_k\cdot\alpha_k\mid k\in\mathbb{N},x_i\in I,\alpha_i\in \SMQ,1\leq i\leq k \}$$ (when $k=0$, we define $x_1\cdot\alpha_1+\cdots+x_k\cdot\alpha_k = 0$).

Assume $\MN$ is a connected marked three manifold.
Let $\pi_1\MN$ be the fundamental groupoid of $\MN$ when $\mathcal{N}\neq\emptyset$, see Definition 3.8 in \cite{wang2023stated}. Then we define $\mathcal{X}\MN$ to be the set of homomorphisms from 
$\pi_1\MN$ to $SL_2(\mathbb{C})$ when $\mathcal{N}\neq\emptyset$, and define $\mathcal{X}\MN$ to be $\mathcal{X}(M)$ when $\mathcal{N}=\emptyset$.
 Then $\XMN$ is an algebraic set.
For any $\rho\in\XMN$, there is an algebra homomorphism $f_{\rho}:\SMN\rightarrow \mathbb{C}$. This offers  
 a one to one correspondence between $\XMN$ and the set of algebra homomorphisms from $\SMN$ to $\mathbb{C}$
\cite{wang2023stated}.
Although we assume $M$ to be connected,
the above discussion can be easily and naturally generalized to general marked three manifolds.

For any $\rho\in \XMN$, we define the {\bf reduced stated skein module} $S_{q^{1/2}}(M,\mathcal{N})_{\rho}$ of $\MN$ with respect to $\rho$ to be 
$$  S_{q^{1/2}}(M,\mathcal{N})/I_{\rho} $$
where $I_{\rho} = (\text{Ker}f_{\rho})\SMQ$.
 When $\mathcal{N}$ is empty,  $S_{q^{1/2}}(M,\mathcal{N})_{\rho}$ coincides with the reduced skein module in equation \eqref{def_reduced}. 

%
%
%
%
%
%
Then we have the parallel statements for reduced stated skein module as in Remark \ref{rem3.6}.
Note that this definition coincides with the definition for  reduced stated $SL_n$-skein module  in Remark 8.17 in \cite{wang2023stated} when $n=2$. 

Then Theorem \ref{th1.1} clearly indicates the following Theorem.

\bt\label{reduced}
Let $(M,\mathcal{N})$ be a compact marked three manifold, and $\rho\in\XMN$.
Then $S_{q^{1/2}}(M,\mathcal{N})_{\rho}$ is finite dimensional, and 
$$\text{dim}_{\mathbb{C}}S_{q^{1/2}}(M,\mathcal{N})_{\rho}  \leq\text{dim}_{\SMN}S_{q^{1/2}}(M,\mathcal{N}).$$

\et

\subsection{Generalization to stated $SL_n$-skein modules}
The author generalized Frobenius map to stated $SL_n$-skein modules, and conjectured the stated $SL_n$-skein modules are finitely generated over the image of the  Frobenius map \cite{wang2023stated}.
The technique used this section can be easily generalized to $SL_n$ case. So we have the following conclusion.

\bt
Conjecture 8.16 in \cite{wang2023stated} holds when the marked three manifold contains at least one marking at every component.
\et

\bt
Let $(M,\mathcal{N})$ be a compact marked three manifold such that every component of $M$ contains at least one marking. Then the reduced
stated $SL_n$-skein module of $(M,\mathcal{N})$, defined in  in Remark 8.17 in \cite{wang2023stated}, is finite dimensional.
\et

\def \MN {(M,\mathcal{N})}
\def \MNP {(M,\mathcal{N}^{'})}
\def \MNO {(M,\mathcal{N}_1)}
\def \MNT {(M,\mathcal{N}_2)}
\def \MNi {(M,\mathcal{N}_i)}
\def \S {S_{q^{1/2}}}

\def \dimQ {\dim(M,Q,q^{1/2})}

\def \dimg {\dim(g,k,q^{1/2})}

\def \vk {\vec{k}}
\def \bZp {\bZ^{\oplus p}}

\def \bN{\mathbb{N}}

\def \bNp{\mathbb{N}^{\oplus p}}
\def \Np{\mathbb{N}^{\oplus p-1}}
\def \vk {\vec{k}}
\def \va {\vec{a}}
\def \vb {\vec{b}}

\section{Dimension of $C(\Si)$ over $\SiN$ when $\Sigma$ is a closed pb surface}

In this section we always assume $\Sigma$ is a connected closed pb surface, and $q^2$ is a primitive $N$-th root of unity with $N$ being odd. The following is the main result of this section.

\bt\label{center}
Let $\Sigma$ be a surface with $p$ punctures and with negative Euler characteristic. Then $C(\Si)$ is freely generated over
$\SiN$ by $N^{p}$ elements. Especially we have $$\dim_{\SiN}C(\Si) = N^{p}.$$
\et

Recall that
for any nonnegative integer $k$,  $P_k$  denotes the surface obtained from $S^{2}$ by removing $k$ points.

\bcr\label{handlebody}
We have $S_q(P_3)$ is freely generated over $S_q(P_3) ^{(N)}$ by $N^3$ elements.
Especially $\dim_{S_q(P_3) ^{(N)}} S_q(P_3)  = N^3$.
\ecr
\bpr
We have $C(S_q(P_3)) = S_q(P_3)$ \cite{le2006colored,przytycki2006skein}. Then Corollary \ref{handlebody} comes from Theorem \ref{center}.
\epr

\begin{remark}

The basis for $C(\Si)$ in
Theorem 3.5 in \cite{frohman2021dimension} can  indicate Theorem \ref{center}. In Remark 5.4 in \cite{frohman2019unicity},
Charles Frohman, Joanna Kania-Bartoszynska, and Thang L{\^e} also claimed the result in Theorem \ref{center}.
 
\end{remark}

In the following of this section, we give a proof for Theorem \ref{center} using quantum trace map.

\subsection{Chekhov-Fock algebra and Quantum trace map}
Let $\Sigma$
be a surface with genus $g$ and $p$ punctures, where $p>0$. 
For an ideal triangulation $\tau$ and a nonzero complex number $\mu$, we have the associated Chekhov-Fock algebra $\mathcal{T}^{\mu}_{\tau}$ see \cite{bonahon2011quantum,bonahon2016representations,bonahon2017representations,liu2009quantum} for more details. Suppose the set of edges of $\tau$ is $\{e_1,\cdots,e_n\}$ where $n = 6 g + 3 p -6$.
 Set $b_{ij}$ to be the number of times an
end of the edge $e_j$
immediately succeeds an end of $e_i$ when going counterclockwise around a
puncture of $\Sigma$, and define $\sigma_{ij} = b_{ij}-b_{ji}\in\{-2, -1, 0, 1, 2\}$.
As an algebra $\mathcal{T}^{\mu}_{\tau}$  is generated by $Y_1^{\pm1},Y_2^{\pm1},\dots,Y_n^{\pm1}$ and subject to the relations: $$Y_i Y_i^{-1} = Y_i^{-1}Y_i = 1,
Y_i Y_j = \mu^{2\sigma_{ij}} Y_j Y_i,$$
where each $Y_i$ is associated to the edge $e_i$.
 For any $Y_{i_1},Y_{i_2},\dots, Y_{i_k}$, we use $[Y_{i_1}Y_{i_2}\dots Y_{i_k}]$ to donate
 $$\mu^{-\sum_{1\leq j<l\leq k}\sigma_{jl}}Y_{i_1}Y_{i_2}\dots Y_{i_k}.$$
For any $\vec{k}=(k_1,\cdots,k_n)\in \mathbb{Z}^{\oplus n},$ we use $Y^{\vec{k}}$
to denote $[Y_{1}^{k_1}Y_{2}^{k_2}\dots Y_{n}^{k_n}]$. Then the set $\{Y^{\vec{k}}\mid \vec{k} \in \mathbb{Z}^{\oplus n}\}$ is a basis of $\Tu$.

For any $\vec{k}=(k_1,\cdots,k_n)\in \mathbb{Z}^{\oplus n}$, we say $\vec{k}$ satisfies the balanced condition if for every ideal triangle $T$ in $\tau$ with three edges labeled by $i_1,i_2,i_3$ we have $k_{i_1}+k_{i_2}+k_{i_3}$ is even.
We use $\Zu$, called the balanced Chekhov-Fock algebra, to denote a subalgebra of $\Tu$ generated by all 
$Y^{\vec{k}}$ with $\vec{k}$ satifying the balanced condition.

 Suppose $v$ is a puncture of $\Sigma$, we use $H_{v}$ to denote $[Y_{i_1}Y_{i_2}\dots Y_{i_{k_v}}] \in \Tu$ 
 where $e_{i_1},e_{i_2},\dots, e_{i_{k_v}}$ are all the edges that connect to the vertex $v$ (note that there maybe some repetitions among $e_{i_1},e_{i_2},\dots, e_{i_{k_v}}$ ). It is easy to show $H_v\in C(\Tu)$.

Set $\mu = q^{1/2}$. Then
there is an algebraic embedding \cite{bonahon2011quantum}
$$Tr_{\tau}^{q}:S_q(\Sigma)\rightarrow \Tu.$$

We set $\nu = \mu^{N^2}$, then there is an algebraic embedding $F:\Tv\rightarrow \Tu$ defined by 
$F(Y^{\vec{k}}) = Y^{N\vec{k}}$.  We also have $\nu^2 = \mu^{2N^2} = q^{N^2} = q^N = \zeta$.
Then we have the following commutative diagram \cite{bonahon2016representations}:
\begin{equation}
\begin{tikzcd}
\Siz  \arrow[r, "\Tz"]
\arrow[d, "\cF"]  
&  \Tv  \arrow[d, "F"] \\
 \Si  \arrow[r, "\Tq"] 
&  \Tu\\
\end{tikzcd}.
\end{equation}
Actually the image of the quantum trace map lies in the balanced Chekhov-Fock algebra. Obviously $F$ restricts to an algebraic embedding $F:\Zv\rightarrow \Zu$. Thus we can replace $\Tv$ (respectively $\Tu$) with $\Zv$ (respectively $\Zu$) in the above diagram, and the new diagram is still commutative.

\subsection{One grading for the balanced Chekhov-Fock algebra}\label{balanced_ch}
We suppose all the punctures of $\Sigma$ are $v_1,\dots,v_p$. And for each $1\leq i\leq p$ we use $H_i$ to denote $H_{v_i}$.

\bl\label{diagonal}
There exist elements $Z_1,\dots,Z_n\in \Zu$ such that (1) for $1\leq i\leq p$ we have $Z_i = H_i$, (2) for any pair $1\leq i<j\leq n$ we have $Z_i Z_j = \mu^{2b_{ij}} Z_j Z_i$ where $b_{ij}\in\mathbb{Z}$, (3)
$\{Z^{\vec{k}}\mid \vec{k} \in \mathbb{Z}^{\oplus n}\}$ is a basis of $\Zu$ where $Z^{\vec{k}}$ is defined in the same way as $Y^{\vec{k}}$.
\el

\bpr
It is an immediate consequence from Lemma 12 in \cite{bonahon2017representations}.
\epr


We suppose $T_N  = \sum_{0\leq t\leq N} \lambda_t x^t$. Note that $\lambda_N  = 1$.

 For any $Z^{\vec{k}}$ with $\vec{k} = (k_1,\cdots,k_n)$, 
we define $\sigma(Z^{\vec{k}}) = (k_1,\cdots, k_p) \in\bZ^{\oplus p}$. For any element 
$x\in \bZ^{\oplus p}$, let $D_x$ be the subvector space of $\Zu$ spanned by all $Z^{\vec{k}}$ with $\sigma(Z^{\vec{k}})=x$. Then we have $\Zu  = \oplus_{x\in \bZ^{\oplus p}} D_x$. Note that this grading is compatible with the algebra structure for $\Zu$, that is, $D_x D_y \subset D_{x+y}$.

For any two different elements $a=(a_1,\cdots,a_p), b=(b_1,\cdots,b_p)\in \bZ^{\oplus p}$, there exists $1\leq i\leq p$ such that $a_t = b_t$ when $1\leq t\leq i-1$ and $a_i\neq b_i$, then we define $a<b$ (respectively $b<a$) if $a_i<b_i$ (respectively $b_i<a_i$). Simply speaking, the linear order "$\leq$" on $\bZ^{\oplus p}$ is lexicographic order.

\def \deg {\text{deg}}

For any nonzero element $l\in \Zu$, we suppose $l = l_{a_1}+\cdots + l_{a_k}$ where $a_i,1\leq i\leq k,$ are $k$ distinct elements in $\bZ^{\oplus p}$ and $l_{a_i}$ is a nonzero element in $ D_{a_i},1\leq i\leq k.$  Then we define 
$\text{deg}(l) = \text{max}\{a_1,\cdots,a_k\}$. For any two nonzero elements $l_1,l_2\in \Zu$,  clearly we have 
$\deg(l_1l_2) = \deg(l_1) + \deg (l_2)$.

For $\Zv$, we define the same $Z_1,\cdots,Z_n\in\Zv$ as $\Zu$ in Lemma \ref{diagonal}. That is, for each $1\leq i\leq n$, suppose $Z_i = Y^{\vec{k_{i}}}\in\Zu$, then similarly we define $Z_i = Y^{\vec{k_{i}}}\in\Zv$. Recall that we have an algebraic embedding $F:\Zv\rightarrow \Zu$, which is a restriction of $F:\Tv\rightarrow \Tu$. Then for each $1\leq i\leq n$, we have $F(Z_i) = F(Y^{\vec{k_{i}}}) = Y^{N\vec{k_{i}}} = (Y^{\vec{k_{i}}})^N = (Z_{i})^N.$

\subsection{Proof of Theorem \ref{center}}
When $p=0$, we have $C(\Si) = \SiN$ \cite{frohman2019unicity}, then it is clearly true for this case.

Then we look at the case when $p>0$. Here we use  all the notations in subsection \ref{balanced_ch}. For each puncture $v_i$, we use $d_i$ to denote the loop going around $v_i$ with vertical framing. From \cite{frohman2019unicity}, we know $C(\Si)$ is generated by $\SiN$ and $\{d_1,\dots,d_p\}$ as a subalgebra of $\Si$. For any $\vec{k}= (k_1,\cdots, k_p)\in\mathbb{N}^{\oplus p}$, we use $d^{\vec{k}}$ to denote
$(d_1)^{k_1}\cdots (d_p)^{k_p}$. Let $\mathcal{C} = \{(k_1,\cdots, k_p)\in\mathbb{N}^{\oplus p}\mid 0\leq k_i\leq N-1, 1\leq i\leq p\}$. Then we want to show $C(\Si)$ is freely generated by $\{d^{\vec{k}}\mid
\vec{k}\in \mathcal{C} \}$ over $\SiN$.

First we want to show $C(\Si)$ is generated by $\{d^{\vec{k}}\mid
\vec{k}\in \mathcal{C} \}$ over $\SiN$. Let $V$ be the $\SiN$-submodule of $C(\Si)$ generated by 
$\{d^{\vec{k}}\mid
\vec{k}\in \mathcal{C} \}$. Then it suffices to show $d^{\vk}\in V$ for all $\vk\in \bNp$. For any $\vk\in\Np$,
$1\leq t\leq p, a\in\bN$, we define $\vk^{t}_{a} = (k_1,\cdots,k_{t-1},a,k_t,\cdots,k_{p-1})\in \bNp$.
For any $\vk\in\Np$ with $0\leq k_i\leq N-1$ for all $ 1\leq i\leq p-1$ and 
$a\in\bN$, we want to show $d^{\vk^1_{a}}\in V$. We already have $d^{\vk^1_{a}}\in V$ when $0\leq a\leq N-1$.
Then suppose we have $d^{\vk^1_{a}}\in V$ when $a\leq m$ where $m\geq N-1$. From the assumption, we have $d^{\vk^1_{m+1-N}}\in V$, then we have 
$$T_N(d_1) d^{\vk^1_{m+1-N}} = \sum_{0\leq t\leq N} \lambda_t d_1^{t}  d^{\vk^1_{m+1-N}} 
= \sum_{0\leq t\leq N} \lambda_t   d^{\vk^1_{t+m+1-N}} \in V. $$
From the assumption, we have $d^{\vk^1_{t+m+1-N}} \in V$ when $0\leq t\leq N-1$, then $d^{\vk^1_{m+1}}\in V$.
From  mathematical induction, we get $d^{\vk^1_{a}}\in V$ when  $\vk\in\Np$ with $0\leq k_i\leq N-1$ for all $ 1\leq i\leq p-1$ and 
$a\in\bN$. Using the same trick as above, we can prove $d^{\vk^2_{a}}\in V$ when  $\vk\in\Np$ with $ 0\leq k_i\leq N-1$ for all $ 2\leq i\leq p-1$ and 
$a\in\bN$. Repeat the above process, eventually we can prove $d^{\vk^p_{a}}\in V$ when  $\vk\in\Np$ and
$a\in\bN$. This actually shows  $d^{\vk}\in V$ for all $\vk\in \bNp$.

Then we try to show the independence of $\{d^{\vec{k}}\mid
\vec{k}\in \mathcal{C} \}$. Suppose $$\sum_{\vec{k}\in \mathcal{C}}\cF(l_{\vec{k}}) d^{\vec{k}} = 0$$
where $l_{\vec{k}}\in \Siz$. Then we want to show $l_{\vec{k}} = 0$ for all $\vec{k}\in \mathcal{C}$. Suppose on the contrary, then $\mathcal{C}_0 = \{\vec{k}\in\mathcal{C}\mid l_{\vec{k}}\neq 0\} \neq \emptyset$.
 We have
\begin{align*}
 0&= \Tq( \sum_{\vec{k}\in \mathcal{C}}\cF(l_{\vec{k}}) d^{\vec{k}})
= \sum_{\vec{k}\in \mathcal{C}_0} \Tq \cF(l_{\vec{k}}) \Tq(d^{\vec{k}})\\ &= 
\sum_{\vec{k}\in \mathcal{C}_0} F \Tz (l_{\vec{k}}) (Z_1 + Z_1^{-1})^{k_1}\cdots (Z_p + Z_p^{-1})^{k_p}\text{ where } \vec{k} = (k_1,\cdots,k_p).
\end{align*}
Since, for each $\vec{k}\in\mathcal{C}_0$, $l_{\vec{k}}\neq 0$, then deg$(F\Tz(l_{\vec{k}}))$ is well-defined. And we have deg$(F\Tz(l_{\vec{k}})) = N x_{\vk}$ for some $x_{\vk}\in \bZp$ because of the definition of $F$.
Then $$\deg(F \Tz (l_{\vec{k}}) (Z_1 + Z_1^{-1})^{k_1}\cdots (Z_p + Z_p^{-1})^{k_p}) = N x_{\vk} + \vk.$$
For any $\va\neq \vb\in \mathcal{C}_0$, we have $N x_{\va} + \va \neq N x_{\vb} + \vb$. Otherwise we have 
$\va -\vb = N (x_{\vb} - x_{\va} )$. But $0\leq a_i, b_i\leq N-1$, thus we get $\va -\vb = N (x_{\vb} - x_{\va} ) = 0$, which contradicts with $\va \neq \vb$.
Thus we have $\deg(\Tq( \sum_{\vec{k}\in \mathcal{C}}\cF(l_{\vec{k}}) d^{\vec{k}}))
=\text{max}\{ N x_{\vk} + \vk\mid \vk\in\mathcal{C}_0 \}$, which contradicts with
$0= \Tq( \sum_{\vec{k}\in \mathcal{C}}\cF(l_{\vec{k}}) d^{\vec{k}})$. Then $l_{\vec{k}} = 0$ for all $\vec{k}\in \mathcal{C}$.



\section{Dimension of  stated skein modules over  Frobenius}

We showed the finiteness of stated skein modules over  Frobenius. In this section, we  focus on estimating the dimension of  stated skein modules over  Frobenius.
For simplicity, we will assume all the three manifolds and pb surfaces, mentioned in this section, are connected.

\def \SqSN {S_{q^{1/2}}(\Sigma)^{(N)}}
\def \fSiN {\widetilde{S_{q^{1/2}}(\Sigma)^{(N)}}}
\def \fSi {\widetilde{\SqS}}

Let $\Sigma$ be a pb surface,  when $q^{1/2}$ is a root of unity of order $N$ with $N$ being odd,
 $\SqSN$ is a commutative domain. 
Recall that
we use $\fSiN$ to denote the field of fractions of $\SqSN$, and use 
$\fSi$ to denote $\SqS\otimes_{\SqSN}\fSiN$. 

 We use $K_{\Sigma,q^{1/2}}$ to denote 
$$\dim_{\fSiN}\fSi,$$
and use $\lambda_{\Sigma,q^{1/2}}$ to denote 
$$\dim_{\SqSN} \SqS.$$
Then clearly we have $K_{\Sigma,q^{1/2}}\leq \lambda_{\Sigma,q^{1/2}}$. 
In this section,
we will give explicit calculation for $K_{\Sigma,q^{1/2}}$, and give a lower bound and an upper bound for $\lambda_{\Sigma,q^{1/2}}$.

\def \wSi {\widetilde{\Si}}
\def \wSiN {\widetilde{\SiN}}


\subsection{When the marking set is empty}\label{sub6.1}

In this subsection, we always assume $q^2$ is a primitive $N$-th root of unity with $N$ being odd. 
From Theorem \ref{skein}, we know
for any compact  three manifold $M$,  $\Sq$ is finitely generated over $\SqN$ where $\zeta = q^{N} = \pm1$. We will give an upper bound for $\dim_{\SqN}\Sq$  in this subsection.

For a closed pb surface $\Sigma$, we know $\Si$ is a domain. Then $\SiN$ is a commutative domain. We use 
 $\widetilde{\SiN}$ to denote the field of fractions of $\SiN$, and use $\wSi$ to denote
$\Si\otimes_{\SiN}\widetilde{\SiN}$. Then $\wSi$ is a vector space over the field $\wSiN$.
Clearly we have $\dim_{\SiN}\Si\geq \dim_{\wSiN}\wSi$.


\bl\label{independence_skein}
Suppse $\Sigma$ is a closed pb surface with genus $g$ and $p$ punctures and the Euler characteristic  $\chi(\Sigma)<0$. Then there exist
$N^{6g-6+3p}$ elements in $S_q(\Sigma)$ which are linearly independent over $S_q(\Sigma)^{(N)}$.
\el
\bpr
Corollary 5.2 in \cite{frohman2021dimension} and Theorem \ref{center}.
\epr

\bl\label{lower_skein}
Suppse $\Sigma$ is a closed pb surface with genus $g$ and $p$ punctures and the Euler characteristic  $\chi(\Sigma)<0$. Then 
$\dim_{\SiN}\Si\geq \dim_{\wSiN}\wSi\geq N^{6g-6+3p}$.
\el
\bpr
Lemma \ref{independence_skein}.
\epr

\bt\label{local_skein} (Corollary 3.10 in \cite{frohman2018structure})
Suppse $\Sigma$ is a closed pb surface with genus $g$ and $p$ punctures and the Euler characteristic  $\chi(\Sigma)<0$. Then 
$ \dim_{\wSiN}\wSi= N^{6g-6+3p} = N^{-3\chi(\Sigma)}$. 
\et
%

 Thang L{\^e} informed us that  Theorem 6.1 in \cite{frohman2021dimension}  indicates the above Theorem.

\bl\label{torus}
Let $T$ be the solid torus. Then we have $S_q(T)$ is a free $S_q(T)^{(N)}$-module   generated by
$N$ elements. Especially, we have $$\dim_{\widetilde{S_q(T)^{(N)}}}\widetilde{ S_q(T)}=\dim_{S_q(T)^{(N)}} S_q(T)= \dim_{S_{\zeta}(T)} S_q(T) = N.$$
\el
\bpr
It is well-known  that $S_q(T) = \mathbb{C}[x]$ and $S_q(T)^{(N)} = \mathbb{C}[T_N(x)]$. Clearly
$S_q(T)$ is  freely generated by $1,x,\cdots,x^{N-1}$  over $S_q(T)^{(N)}$.

\epr

\begin{remark}
There are only four closed pb surfaces with non-negative Euler characteristic: closed torus $T,P_0,P_1,P_2$.
Proposition 5.5 in \cite{frohman2016frobenius} shows $\dim_{\widetilde{S_q(T)^{(N)}}} \widetilde{S_q(T)} = N^2$.
Obviously we have $\dim_{\widetilde{S_q(P_i)^{(N)}}} \widetilde{S_q(P_i)} = 1$ when $i=0,1$.
Lemma \ref{torus} shows $\dim_{\widetilde{S_q(P_2)^{(N)}}} \widetilde{S_q(P_2)} = N$.
In \cite{frohman2016frobenius}, Frohman and Abdiel calculated 
$\dim_{\widetilde{S_q(\Sigma)^{(N)}}} \widetilde{S_q(\Sigma)}$ when 
$\Sigma$ is the closed torus, once punctured torus, $P_2,P_3$, which coincide with Theorem \ref{local_skein} and Lemma \ref{torus}.
\end{remark}


\bt(First paragraph in the proof of Theorem 6.2 in  \cite{przytycki1997skein})\label{finiteorder}
Suppse $\Sigma$ is a closed pb surface with genus $g$ and $p$ punctures, where $p\geq 1$. Then there exist closed curves $\alpha_1,\dots,\alpha_{k}$ such that $\Si$ is orderly finitely generated by $\alpha_1,\dots,\alpha_{k}$, where $k=2^{2g + p -1} -1$.
\et

\bt\label{skein_l_u}
Suppse $\Sigma$ is a closed pb surface with genus $g$ and $p$ punctures, where $p\geq 1$, and the Euler characteristic  $\chi(\Sigma)<0$.
Then $$N^{6g-6+3p}\leq \dim_{\SiN}\Si\leq N^{2^{2g + p -1}-1}.$$
\et
\bpr
 From Theorem \ref{finiteorder},
we know there exist $k=2^{2g + p -1} -1$ disjoint solid tori embedded in $\Sigma\times[0,1]$ such that the embedding $f$ from the union of these $k$ solid tori, which is denoted as $M^{'}$, to $\Sigma\times[0,1]$ induces a surjective linear map
$\fs:S_q(M^{'})\rightarrow S_q(\Sigma\times[0,1])$. Then from Lemma \ref{keylemma}, we have 
$$\dim_{\SiN}\Si = \dim_{S_q(\Sigma\times [0,1])^{(N)}} S_q(\Sigma\times [0,1])\leq 
\dim_{S_q(M^{'})^{(N)}} S_q(M^{'}) = N^{k}$$
where the last equality comes from Lemma \ref{torus}.
Then Lemma \ref{lower_skein} completes the proof.
\epr

\bt
Suppse $\Sigma$ is a totally closed pb surface with genus $g$ and the Euler characteristic  $\chi(\Sigma)<0$.
Then $$N^{6g-6}\leq \dim_{\SiN}\Si\leq N^{2^{2g}-1}.$$
\et
\bpr

Let $\Sigma^{'}$ be a pb surface obtained from $\Sigma$ by removing one point in $\Sigma$. Then the embedding
$l:\Sigma^{'}\rightarrow \Sigma$ induces a surjective algebraic map $l_{\sharp}:S_q(\Sigma^{'})\rightarrow S_q(\Sigma)$. Then Lemma \ref{keylemma} and Theorem \ref{skein_l_u} indicate
$$\dim_{\SiN}\Si\leq\dim_{S_q(\Sigma^{'})^{(N)}}S_q(\Sigma^{'})\leq N^{2^{2g}-1}.$$
Then Lemma \ref{lower_skein} completes the proof.
\epr

\bt\label{thm6.7}
Let $M$ be a  compact three manifold with genus $g$,
%
%
%
%
we have $$\dim_{S_{\zeta}(M)} \Sq \leq N^{2^g-1}.$$ Especially the dimension of any reduced skein module of $M$
is not more than $N^{2^g-1}$.
\et
\bpr
Equation \eqref{eq_handle}, Lemma \ref{torus}  and Theorems \ref{reduced_dim}, \ref{skein_l_u}.
\epr

\bcr
Let $M$ be a  compact   three manifold with genus $g$. If $\Sz$ is finitely generated over $\Szp$ where $\zeta = q^N = \pm1$, then 
$$\text{dim}_{\Sqp}\Sq \leq N^{2^g-1} (\text{dim}_{\Szp}\Sz).$$
\ecr
\bpr
Theorems \ref{thm3.9} and \ref{thm6.7}.
\epr

In the remaining of this subsection, we will   focus on proving $\cF:S_{\zeta}(S^1\times S^2) \rightarrow S_q(S^1\times S^2)$ is an isomorphism, which indicates
 $S_q(S^1\times S^2)$ is freely  generated over $S_{\zeta}(S^1\times S^2)$ by the empty skein.

Recall that $T_n(x)$ and $S_n(x)$ are Chebyshev polynomials of the first kind and of the second kind respectively.
Now we define a new sequence of  polynomials $A_n(x)$ by setting $A_n(x) = S_n(x),n=1,2,$ and $A_n(x) = S_n(x) + A_{n-2}(x)$ when $n>2$. Suppose $D$ is an embedded disk in $S^2$, then $S^1\times D$ is an embedded solid torus in $S^1\times S^2$. We use $P$ to denote the origin of the disk $D$, and use $x$ to denote the skein in $S^1\times D$ represented by the closed line $\{P\}\times S^1$ with vertical framing. It is well-known that $S_q(S^1\times D)$
is actually $\mathbb{C}[x]$. 

For any polynomial $f(x)\in \mathbb{C}[x]$, we can regard $f(x)$ as an element in 
$S_q(S^1\times S^2)$ by the embedding from $S^1\times D$ to $S^1\times S^2$. Then it is obvious that $1,x,x^2,\dots$ span $S_q(S^1\times S^2)$, and $1,T_N(x),T_N(x)^2,\dots$ span $S_q(S^1\times S^2)^{(N)}$.

\bt
$\cF:S_{\zeta}(S^1\times S^2) \rightarrow S_q(S^1\times S^2)$ is an isomorphism. Especially, we have
$S_q(S^1\times S^2)$ is freely generated over $S_{\zeta}(S^1\times S^2)$ by the empty skein, and $S_q(S^1\times S^2)$ has a commutative algebra structure.
\et
\bpr
If $N$ is 1, it is obvious.

\def \Ns {\mathbb{N}^{*}}

Suppose $N\geq 3$. 
From
\cite{hoste1995kauffman}, we know $S_{\zeta}(S^1\times S^2)$ has a basis $$\{1,A_1(x),\cdots,A_k(x),\cdots\}.$$
Thus $\{1,T_1(x),\cdots,T_k(x),\cdots\}$ is also a basis for $S_{\zeta}(S^1\times S^2)$. For any $k\in\mathbb{N}^{*}$, we have $\cF(T_k(x)) = T_{kN}(x)$. Thus it suffices to show 
$\{1,T_{N}(x),\cdots,T_{kN}(x),\cdots\}$ is a basis for $S_q(S^1\times S^2)$.


For each positive integer $i$, we set $e_i$ to be $A_i(x)\in S_q(S^1\times S^2)$.
From
\cite{hoste1995kauffman}, we know 
$$S_q(S^1\times S^2) = \mathbb{C}\emptyset \oplus_{1\leq i<+\infty, N\mid i+2} \mathbb{C}e_i,$$
and $e_i = 0$ for $1\leq i<+ \infty, N\nmid i+2$. Then 
\begin{equation}\label{eq_new}
\{1\}\cup \{e_i\mid \exists k\in\Ns \text{ such that } i+2 = kN\}
\end{equation}
is a basis for $S_q(S^1\times S^2)$.

Suppose $i+2 = kN$ for some positive integer $k$. 
We have 
\begin{align*}
T_{i+2}(x) &= S_{i+2}(x) - S_{i}(x) = A_{i+2}(x) - A_i(x) -(A_i(x) - A_{i-2}(x))\\
&= A_{i+2}(x) + A_{i-2}(x) - 2A_i(x) = e_{i+2} + e_{i-2} - 2e_i.
\end{align*}
Since $N\mid i+2$, $N$ is odd and $N$ is not 1, we have $N\nmid i+4$ and $N\nmid i$. Then 
$e_{i+2} = e_{i-2} = 0$. Thus we have $T_{i+2}(x) = -2e_i$, then $e_i = -\frac{1}{2} T_{kN}(x)\in  S_q(S^1\times S^2).$  Note that there is a special case where $N=3,k=1$, but using the same technique as above we can still get $e_i = -\frac{1}{2} T_{kN}(x)\in  S_q(S^1\times S^2).$
Then basis in equation \eqref{eq_new} indicates $\{1,T_{N}(x),\cdots,T_{kN}(x),\cdots\}$ is a basis for $S_q(S^1\times S^2)$.
\epr

\begin{corollary}\label{S1S2}
The reduced skein module for $S^1\times S^2$ always has dimension one.
\end{corollary}

\begin{remark}
Charles Frohman, Joanna Kania-Bartoszynska, and Thang L{\^e} prove the reduced skein module of any closed three manifold with respect to any non-central character is one dimensional \cite{talk}.

Corollary \ref{S1S2} indicates the reduced skein module for $S^1\times S^2$ with respect to any central character is still one dimensional.
\end{remark}

%
%

\def \Oq {O_q(SL_2)}
\def \OqS {\widetilde{O_q(SL_2)}}
\def \Kq {K_{q^{1/2}}}
\def \lmq {\lambda_{q^{1/2}}}

\subsection{On $\dim_{S_{q^{1/2}}(D_2)^{(N)}}S_{q^{1/2}}(D_2)$}
In the remaining of this section, we will always assume $q^{1/2}$ is a root of unity of order $N$ with $N$ being odd.
Recall that
$D_2$ is the  bigon, and  $\lambda_{q^{1/2}}$  denotes
$\text{dim}_{S_{q^{1/2}}(D_2)^{(N)}}S_{q^{1/2}}(D_2)$.
From \cite{costantino2022stated1}, we know $S_{q^{1/2}}(D_2)$ is just $O_q(SL_2)$, where
$\Oq$ is generated by $a,b,c,d$ and subject to the following relations:
\begin{align*}
ca &= q^2ac, db = q^2bd, ba = q^2ab, dc=q^2cd,\\
bc&=cb, ad - q^{-2}bc =1, da-q^2 cb =1.
\end{align*}
Let $A_{q}$ be the the subalgebra of $\Oq$ generated by $a^N,b^{N}, c^N,d^N.$ With the identification between 
$S_{q^{1/2}}(D_2)$ and $O_q(SL_2)$, $S_{q^{1/2}}(D_2)^{(N)}$ is just $A_q$. Thus we have 
 $\lambda_{q^{1/2}}=\text{dim}_{A_q}\Oq$.
 Since $\Oq$ is a domain, then 
$A_q$ is a commuative domain. We use $\widetilde{A_q}$ to denote the field of fractions of $A_q$, and 
use $\OqS$ to denote 
$ \Oq\otimes_{A_q} \widetilde{A_q}$. We use $\Kq$ to denote
$$\dim_{\widetilde{A_q}}\OqS = \dim_{\widetilde{S_{q^{1/2}}(D_2)^{(N)}}}\widetilde{S_{q^{1/2}}(D_2)}.$$ Clearly we have $K_{q^{1/2}}\leq \lmq$.

\begin{remark}

$K_{q^{1/2}}$ is actually the rank of $\Oq$ over $O(SL_2)$. In \cite{brown2012lectures} chapter 3,  there is a general calculation for this rank for all semisimple Lie groups (our case is $SL_2$). In this subsectoin, we will give an elementary way to calculate $K_{q^{1/2}}$, and give  a clear basis for $\widetilde{\Oq}$ over $\widetilde{A_q}$.

\end{remark}

Since the map from $\Oq$ to $\OqS$, given by $x\mapsto x\otimes 1$, for $x\in\Oq$, is injective, we can regard any element in $\Oq$ as an element in $\OqS$ by this embedding.

\bl\label{local}
Suppose $\alpha_1,\dots,\alpha_n\in\Oq$  are linearly independent over $A_q$. When we regard all $\alpha_i$ as elements in $\OqS$, $\alpha_1,\dots,\alpha_n$ are are linearly independent over $\widetilde{A_q}$. Especially we have $\lmq\geq\Kq\geq n$.
\el

\def \bN {\mathbb{N}}

We define $\Lambda = \{(k_1,k_2,k_3,k_4) \mid k_1,k_2,k_3,k_4\in\mathbb{N}, k_1k_2 = 0\}$. For any 
$$\vec{k}= (k_1,k_2,k_3,k_4)\in \bN\times \bN\times \bN\times \bN,$$ we define
$O_{\vec{k}}$ to be $a^{k_1}d^{k_2}b^{k_3}c^{k_4}\in\Oq.$
\bl(\cite{gavarini2007pbw,wang2023stated})\label{basis}
The set $\{O_{\vec{k}}|\vec{k} \in\Lambda\}$ is a basis for 
$\Oq$, and the set 
$\{O_{N\vec{k}}|\vec{k} \in\Lambda\}$ is a basis for 
$A_q$.
\el

For any positive integer $t$,
define a subset of $\Oq$ as $$E_t = \{k_t b^tc^t +\cdots +k_1 bc + 1\mid  k_t = q^s\text{ for some integer } s \text\}.$$ 
We define $E_0 = \{1\}.$ Then for any $l,t\in\mathbb{N}$, we have $E_l E_t\subset E_{l+t}$.

\bl\label{commu}
For any $t\in\mathbb{N}$, define $d E_t
=\{df\mid f\in E_t\},  E_td = \{fd\mid f\in E_t\}.$ Then we have 
we have $d E_t= E_td.$ 
Similarly, we also have $aE_t = E_ta$.
\el
\bpr
Note that $d(bc)=q^4 (bc) d$, which obviously indicates the Lemma.
\epr

\bl\label{a_d}
For any $t\in\mathbb{N}$,
we have $a^t d^t\in E_t$.
\el
\bpr
We prove this Lemma by using mathematical induction on $t$. Note that $ad = q^{-2}bc + 1$, then the Lemma is obviously true when $t = 0,1$.
We suppose $a^td^t = f\in E_t$, then we want to show $a^{t+1}d^{t+1}\in E_{t+1}$.
We have $$a^{t+1} d^{t+1} = afd = ad g
= (q^{-2}bc + 1) g\in E_{t+1}$$ where $g\in E_t$.
\epr

\def \deg {\text{deg}}

We use $\leq_{O}$ to denote the lexicographic order on $\bN\times \bN\times \bN\times \bN$. 
Then for any nonzero element $x\in\Oq$, we have 
$x =\sum_{\vec{k}\in\Lambda_x} l_{\vec{k}} O_{\vec{k}} $
where $\Lambda_x$ is a finite subset of $\Lambda$ and 
$l_{\vec{k}} \in\mathbb{C}^{*}$ for all $\vec{k}\in \Lambda_x$.
Then we define $\deg(x) =\text{max}(\Lambda_x)$
where $\text{max}(\Lambda_x)$ is the maximal element in $\Lambda_x$ under the linear order $\leq_{O}$.
Then we have the following two Lemmas.
\bl\label{degree}
Let $x_1,\cdots,x_k$ be $k$ ($k>0$) nonzero elements in $\Oq$. Suppose
$\deg(x_1),\cdots,\deg(x_k)$ are $k$ distinct elements in $\Lambda$, then
$x_1+\cdots +x_k\neq 0.$
\el
\bpr
It is obvious.
\epr

We define a function $\varphi :\bN\times \bN\times \bN\times \bN\rightarrow \Lambda$ by
$$\varphi(\vec{k}) = \left\{ 
    \begin{aligned}
    &(0,0,k_3+k_1,k_4+k_2) & & k_1 = k_2 \cr 
    &(k_1 - k_2,0,k_3+k_2,k_4+k_2) & & k_1 > k_2 \cr 
    &(0,k_2 - k_1,k_3+k_1,k_4+k_1) & & k_1 < k_2
    \end{aligned}
\right.$$
where $\vec{k} = (k_1,k_2,k_3,k_4)\in \bN\times \bN\times \bN\times \bN.$

\bl\label{Odegree}
For any $\vec{k}\in \bN\times \bN\times \bN\times \bN$, we have 
$\deg(O_{\vec{k}}) = \varphi(\vec{k})$.
\el
\bpr
Let $\vec{k} = (k_1,k_2,k_3,k_4)\in \bN\times \bN\times \bN \times \bN$.

(1) Suppose $k_1 = k_2$. From Lemma \ref{a_d}  we have 
$$O_{\vec{k}} = a^{k_1} d^{k_1} b^{k_3} c^{k_4} = f b^{k_3} c^{k_4}$$
where $f\in E_{k_1}$. Thus $\deg(O_{\vec{k}}) = \deg(f b^{k_3} c^{k_4})
= (0,0, k_3+ k_1,k_4+k_1) = \varphi(\vec{k})$.

(2) Suppose $k_1 > k_2$. From Lemma \ref{a_d}  we have 
$$O_{\vec{k}} = a^{k_1} d^{k_2} b^{k_3} c^{k_4} = a^{k_1-k_2} h b^{k_3} c^{k_4}$$
where $h\in E_{k_2}$. Thus $\deg(O_{\vec{k}}) = \deg(a^{k_1-k_2} h b^{k_3} c^{k_4})
= (k_1-k_2,0, k_3+ k_2,k_4+k_2) = \varphi(\vec{k})$.

(3) Suppose $k_1 < k_2$. From Lemma \ref{a_d}  we have 
$$O_{\vec{k}} = a^{k_1} d^{k_2} b^{k_3} c^{k_4} = x d^{k_2-k_1} b^{k_3} c^{k_4}$$
where $x\in E_{k_1}$. 
From Lemma \ref{commu}, we know there exists $y\in E_{k_1}$ such that 
$x d^{k_2-k_1} =  d^{k_2-k_1} y$. Then 
 $\deg(O_{\vec{k}}) = \deg(d^{k_2-k_1} y b^{k_3} c^{k_4})
= (0,k_2-k_1, k_3+ k_1,k_4+k_1) = \varphi(\vec{k})$.
\epr

There is an obvious partition for $\Lambda$, defined by $\Lambda = \Lambda_0\cup \Lambda_1\cup\Lambda_2$, where 
\begin{align*}
\Lambda_0 &= \{(0,0,k_3,k_4)\mid k_3,k_4\in\mathbb{N} \}\\
\Lambda_1 &= \{(k_1,0,k_3,k_4)\mid k_1\in\mathbb{N}^{*}, k_3,k_4\in\mathbb{N} \}\\
\Lambda_2& = \{(0,k_2,k_3,k_4)\mid k_2\in\mathbb{N}^{*}, k_3,k_4\in\mathbb{N} \}.
\end{align*}

Then for any $\vec{k}\in \bN\times \bN\times \bN\times \bN$, we have
$\varphi(\vec{k})\in\Lambda_0$ if and only if $k_1 = k_2$, 
$\varphi(\vec{k})\in\Lambda_1$ if and only if $k_1 > k_2$,
$\varphi(\vec{k})\in\Lambda_2$ if and only if $k_1 < k_2$.
\bl\label{ob}
For any $\vec{u},\vec{v}\in\Lambda$, if $u_1 - u_2 = v_1 - v_2$, then we have $u_1 = v_1$ and $u_2 = v_2$.
\el
\bpr
Suppose $\vec{u}\in \Lambda_{i},\vec{v}\in\Lambda_{j}$ where $i,j = 0,1,2$. Note that if $i\neq j$, we cannot have
$u_1 - u_2 = v_1 - v_2$,
thus $i=j.$ Then trivially, $u_1 - u_2 = v_1 - v_2$ indicates $u_1 = v_1$ and $u_2 = v_2$.
\epr

\def \cA{\mathcal{D}}

Define a subset $\mathcal{D}$ of $\Lambda$ by
$$\mathcal{D} = \{(0,k_2,k_3,k_4)\in\Lambda\mid 0\leq k_2, k_3,k_4\leq N-1 \}.$$
Define a map $\psi:\Lambda \times \mathcal{D}\rightarrow \Lambda$ by
$\psi(\vec{u},\vec{v}) = \varphi(N\vec{u} + \vec{v})$
where $\vec{u}\in\Lambda,\vec{v}\in \cA$. Then we have the following Lemma.

\bl\label{injective}
The above  map $\psi$ is injective.
\el
\bpr
Assume $\psi(\vec{u},\vec{v}) = \psi(\vec{m},\vec{n})$, then we want to show $\vec{u} =\vec{m}$ and 
$\vec{v} =\vec{n}$.

(1) Suppose $\psi(\vec{u},\vec{v}) = \psi(\vec{m},\vec{n}) \in\Lambda_0$. 
Then we have $Nu_1+v_1 = Nu_2+v_2$ and $Nm_1+n_1 = Nm_2+n_2$.
Since $\vec{v}\in\cA$, we have $v_1 = 0$, then $N(u_1 - u_2)= v_2$, that is $N\mid  v_2$. We also have 
$0\leq v_2\leq N-1$, then $v_2 = 0$. Then we get $u_1 = u_2$. Since $\vec{u}\in\Lambda$, we have 
$u_1 = u_2 = 0$. Similar, we can get $m_1=m_2=n_1=n_2=0.$

Then we have
\begin{align*}
&\psi(\vec{u},\vec{v}) =(0,0,Nu_3 + v_3, Nu_4 + v_4 ) 
= \psi(\vec{m},\vec{n}) =(0,0,Nm_3 + n_3 , Nm_4 + n_4 ).
\end{align*}
For each $i=3,4$, we have $Nu_i + v_i = Nm_i + n_i$. Since $0\leq v_i,n_i\leq N-1$, we have $u_i = m_i,v_i =n_i.$
Thus we get $\vec{u} = \vec{m}$ and $\vec{v} = \vec{n}$.

(2) Suppose $\psi(\vec{u},\vec{v}) = \psi(\vec{m},\vec{n}) \in\Lambda_1$. 
Then we have $Nu_1+v_1 = Nu_1 > Nu_2+v_2$ and $Nm_1+n_1 =Nm_1> Nm_2+n_2$.
 From the defintion of $\psi$, we know
\begin{align*}
\psi(\vec{u},\vec{v})
= &(N(u_1 - u_2)-v_2, 0 , Nu_3 + v_3 + Nu_2+v_2,\\
 &Nu_4 + v_4 +Nu_2+v_2)\\
\psi(\vec{m},\vec{n})
 =& (N(m_1 - m_2)-n_2, 0 , Nm_3 + n_3 +  Nm_2+n_2,\\& Nm_4 + n_4 + Nm_2+n_2).
\end{align*}
From $\psi(\vec{u},\vec{v}) = \psi(\vec{m},\vec{n})$, we get $N(u_1 - u_2)-v_2 = N(m_1 - m_2)-n_2$, then 
$v_2 = n_2$ and $u_1 - u_2 = m_1 - m_2$ since $0\leq v_2,n_2\leq N-1$.
Then Lemma \ref{ob} indicates $u_1 = m_1$ and $u_2 = m_2$. 
From $\psi(\vec{u},\vec{v}) = \psi(\vec{m},\vec{n})$, we can also get
$Nu_i + v_i = Nm_i + n_i$
 for each $i = 3,4$. Similarly as above we can have  $u_i = m_i,v_i =n_i$ for each $i = 3,4$. Thus 
we get $\vec{u} = \vec{m}$ and $\vec{v} = \vec{n}$.

(3) Suppose $\psi(\vec{u},\vec{v}) = \psi(\vec{m},\vec{n}) \in\Lambda_2$. The proving technique is similar with case (2).

\epr

\def \vk {\vec{k}}
\def \vv {\vec{v}}

\bl\label{free1}
We have $\{O_{\vec{k}}\mid \vec{k}\in\mathcal{D} \}$ is linearly independent over $A_q$.
\el
\begin{proof}
It surfices to show 
$$\sum_{\vec{k}\in\cA^{'}}  \alpha_{\vec{k}} O_{\vec{k}} \neq 0$$
where $\cA^{'}$ is a nonempty finite subset of $\cA$ and $ \alpha_{\vec{k}}$ is a nonzero element in $A_q$ for each $\vec{k}\in\cA^{'}$.

For each $\vec{k}\in\cA^{'}$, from Lemma \ref{basis} we can  suppose $\alpha_{\vk} = \sum_{\vv\in\Lambda_{\vk}} l_{\vk,\vv} O_{N\vv}$ where $\Lambda_{\vk}$ is a nonempty finite subset of $\Lambda$ and 
$l_{\vk,\vv}\neq 0$ for each $\vv\in \Lambda_{\vk}$. Then we have
$$\sum_{\vec{k}\in\cA^{'}}  \alpha_{\vec{k}} O_{\vec{k}}
= \sum_{\vec{k}\in\cA^{'}}  \sum_{\vv\in\Lambda_{\vk}} l_{\vk,\vv} O_{N\vv} O_{\vec{k}}= \sum_{\vec{k}\in\cA^{'}}  \sum_{\vv\in\Lambda_{\vk}} l_{\vk,\vv} O_{N\vv +\vk}.$$
From Lemma \ref{Odegree}, we know $\deg(O_{N\vv +\vk}) = \varphi(N\vv +\vk) =\psi(\vv,\vk)$.
From Lemma \ref{injective}, we know $\deg(O_{N\vv +\vk}), \vec{k}\in\cA^{'}, \vv\in\Lambda_{\vk}$, are distinct.
Then from Lemma \ref{degree}, we get $\sum_{\vec{k}\in\cA^{'}}  \alpha_{\vec{k}} O_{\vec{k}}=\sum_{\vec{k}\in\cA^{'}}\sum_{\vv\in\Lambda_{\vk}} l_{\vk,\vv} O_{N\vv +\vk} \neq 0$.
\end{proof}

\bcr
We have $\lmq\geq\Kq\geq N^{3}$.
\ecr
\bpr
Lemmas \ref{local}, \ref{free1} and the fact that $|\cA| = N^{3}$.
\epr

\def \AqS {\widetilde{A_q}}

\begin{theorem}\label{dimen}
When
we regard $\{O_{\vec{k}}\mid \vec{k}\in\mathcal{D} \}$ as a subset of $\OqS$, 
it is a basis of $\OqS$ over $\widetilde{A_q}$.
Especially
we have $\Kq = N^{3}$.
\end{theorem}
\bpr
 Lemmas 
 \ref{local} and \ref{free1} indicate $\{O_{\vec{k}}\mid \vec{k}\in\mathcal{D} \}$ are linearly independent over
$\widetilde{A_q}$. Let $V$ be the linear span of $\{O_{\vec{k}}\mid \vec{k}\in\mathcal{D} \}$ over $\AqS$.
 Then it suffices to show $ V = \OqS$.

 For any $k_1,k_2,k_3\in\mathbb{N},$ clearly we have 
$d^{k_1} b^{k_2} c^{k_3}\in V$. Then we want to show for any $k_1,k_2,k_3,k_4\in\mathbb{N}$ we have 
$a^{k_1} d^{k_2} b^{k_3} c^{k_4}\in V$. Suppose $k_1 = uN + v$ where $u,v\in\mathbb{N}$ and $0\leq v\leq N-1$.
Then we have 
\begin{align*}
&d^N( a^{k_1} d^{k_2} b^{k_3} c^{k_4}) = a^{uN} a^{v} d^N d^{k_2} b^{k_3} c^{k_4} =
a^{uN} a^{v} d^{v} d^{N-v} d^{k_2} b^{k_3} c^{k_4}\\
=& a^{uN} f d^{N +k_2-v}  b^{k_3} c^{k_4} = a^{uN}  d^{N +k_2-v} h b^{k_3} c^{k_4} \in V.
\end{align*}
where $f,h\in E_v$. Since $d^{N}$ is invertible, we have $a^{k_1} d^{k_2} b^{k_3} c^{k_4}\in V$.
\epr

\begin{remark}
During the proof of Theorem \ref{dimen}, we only used the invertibility of $d^{N}$. The same technique can show 
$\{a^{k_1} b^{k_2} c^{k_3}\mid 0\leq k_1,k_2,k_3\leq N-1\}$ span $\OqS$ by using the invertibility of $a^{N}$.
The reduced stated skein module of $D_2$ is defined by  a matrix $\begin{pmatrix}
x_a & x_b\\
x_c & x_d
\end{pmatrix} \in SL_2(\mathbb{C})$. We regard $\begin{pmatrix}
x_a & x_b\\
x_c & x_d
\end{pmatrix} \in SL_2(\mathbb{C})$ as an algebra homomorphism from $A_q$ to $\mathbb{C}$ such that 
it maps $\begin{pmatrix}
a^{N} & b^{N}\\
c^{N} & d^{N}
\end{pmatrix} $to $\begin{pmatrix}
x_a & x_b\\
x_c & x_d
\end{pmatrix} $.
Then the above discussion shows $\{a^{k_1} b^{k_2} c^{k_3}\mid 0\leq k_1,k_2,k_3\leq N-1\}$ (respectively $\{d^{k_1} b^{k_2} c^{k_3}\mid 0\leq k_1,k_2,k_3\leq N-1\}$) spans the corresponding reduced stated skein module of $D_2$  when $x_a \neq 0$ (respectively $x_d \neq 0$).

 We will try to calculate the dimension for reduced stated skein module in future work.
\end{remark}

\def \cB {\mathcal{B}}

Define $$\mathcal{B} =  \{(N-k_1, 0, k_2, k_3)\mid 1\leq k_1\leq N-1, 0\leq k_2,k_3\leq N-1,k_2< k_1\text{ or } k_3< k_1\}\subset \Lambda$$
\bl\label{upper_O}
We have $\Oq$ is linearly spanned by $\{O_{\vk}\mid \vk\in\cA\cup\cB \}$ over $A_q$.
Especially $\lmq\leq |\cA\cup\cB| = 2N^3 - \frac{N(N+1)(2N+1)}{6}$.
\el
\bpr
Let $U$ be the linear span of $\{O_{\vk}\mid \vk\in\cA\cup\cB \}$ over $A_q$. 
Then it suffices to  show 
$a^{N-k_1}b^{k_2}c^{k_3}\in U$ for any $1\leq k_1\leq N-1, 0\leq k_2,k_3\leq N-1.$
From the definition of $\cB$, we know $a^{N-k_1}b^{k_2}c^{k_3}\in U$ if 
$k_2< k_1\text{ or } k_3< k_1$. 
Then we suppose $k_2\geq k_1$ and $k_3\geq k_1$,
and we use $k$ to denote $k_2 -k_1 + k_3-k_1$.
 Then we will prove this case by using mathematical induction on  $k$.
From Lemma \ref{a_d}, we have 
\begin{align*}
a^N b^{k_1} = a^{N-k_1}(l_{k_1}b^{k_1}c^{k_1}+
l_{k_1-1}b^{k_1-1}c^{k_1-1} +\cdots+1)
\end{align*}
where $l_{k_1}$ is a nonzero complex number. Then we have
\begin{align}\label{eee}
a^{N-k_1}b^{k_1}c^{k_1} = s_{k_1} a^N b^{k_1} + s_{k_1 - 1}a^{N-k_1}b^{k_1-1}c^{k_1-1} +\cdots+s_0a^{N-k_1}\in U.
\end{align}
Equation \eqref{eee} indicates $a^{N-k_1}b^{k_2}c^{k_3}\in U$ when $k=0$. Then suppose 
$a^{N-k_1}b^{k_2}c^{k_3}\in U$ when $k=k_2 -k_1 + k_3-k_1\leq m$ ($m\in\mathbb{N}$), then we look at the case when $k= m+1$.
We left multiply equation
\eqref{eee} on both sides by $b^{k_2-k_1}c^{k_3-k_1}$, we get 
$$a^{N-k_1}b^{k_2}c^{k_3} = s_{k_1} a^N  b^{k_2}c^{k_3-k_1} + s_{k_1 - 1}a^{N-k_1}b^{k_2-1}c^{k_3-1} +\cdots+s_0a^{N-k_1} b^{k_2-k_1}c^{k_3-k_1}\in U.$$
\epr

\subsection{When $(M, \mathcal{N})$ is the thickening of a pb surface}


\def \bZp {\bZ^{\oplus p}}

\def \bN{\mathbb{N}}

\def \bNp{\mathbb{N}^{\oplus p}}
\def \Np{\mathbb{N}^{\oplus p-1}}
\def \vk {\vec{k}}
\def \va {\vec{a}}
\def \vb {\vec{b}}

Let $\Sigma$ be a  pb surface. 
Recall that
$r(\Sigma) = -\chi(\Sigma)+\sharp\partial\Sigma$, where $\chi(\Sigma)$ is the Euler characteristic of $\Sigma$ and $\sharp\partial\Sigma$ is the number of boundry components of $\Sigma$. Note that when $\partial\Sigma=\emptyset$, we have $r(\Sigma) = -\chi(\Sigma) = 2g-2+p$ where $g$ is the genus and $p$ is the number of  punctures.
Theorem 7.13 in \cite{le2021stated} indicates the following Lemma.

\bl\label{Saturate}
Let $\Sigma$ be a  pb surface with nonempty boundary.
Then we have a linear isomorphism
$$F : S_{q^{1/2}}(D_2)^{\otimes r(\Sigma)}\rightarrow S_{q^{1/2}}(\Sigma)$$
such that $F$ restricts to an algebraic isomorphism $$F : (S_{q^{1/2}}(D_2)^{(N)})^{\otimes r(\Sigma)}\rightarrow (S_{q^{1/2}}(\Sigma))^{(N)},$$ and $F$ preserves module structures, that is 
$$F(\alpha\cdot\beta) = F(\alpha)\cdot F(\beta)$$ for any $\alpha\in (S_{q^{1/2}}(D_2)^{(N)})^{\otimes r(\Sigma)}, \beta\in S_{q^{1/2}}(D_2)^{\otimes r(\Sigma)}$.
Especially we have 
\begin{align*}
\dim_{S_{q^{1/2}}(\Sigma)^{(N)}}S_{q^{1/2}}(\Sigma) = \dim_{(S_{q^{1/2}}(D_2)^{(N)})^{\otimes r(\Sigma)}} S_{q^{1/2}}(D_2)^{\otimes r(\Sigma)}
=\dim_{(A_q)^{\otimes r(\Sigma)}} \Oq^{\otimes r(\Sigma)}.
\end{align*}
\el

\bl\label{independence}
For any positive integer $k$, $\Oq^{\otimes k}$ contains $N^{3k}$ elements which are linearly independent over 
$(A_q)^{\otimes k}.$
\el
\bpr
From Lemma \ref{free1}, we know there exist $N^3$ elements $x_1,\cdots,x_{N^3}\in\Oq$ which are linearly independent over $A_q$. Then we have $\oplus_{1\leq i\leq N^3} A_q x_i\subset \Oq$.
Thus we have 
\begin{align*}
&(\oplus_{1\leq i\leq N^3} A_q x_i)\otimes\cdots\otimes (\oplus_{1\leq i\leq N^3} A_q x_i)
=\oplus_{1\leq i_1,\cdots,i_{k}\leq N^3} (A_q x_{i_1} \otimes\cdots\otimes  A_q x_{i_{k}})\\
=&\oplus_{1\leq i_1,\cdots,i_{k}\leq N^3} ((A_q)^{\otimes k} x_{i_1} \otimes\cdots\otimes   x_{i_{k}})\subset \Oq^{\otimes k}.
\end{align*}
This shows that $\{x_{i_1} \otimes\cdots\otimes   x_{i_{k}}\mid 1\leq i_1,\cdots,i_{k}\leq N^3 \}$
are linearly independent over $(A_q)^{\otimes k}$.
\epr

\def \SqS {S_{q^{1/2}}(\Sigma)}
\def \SqSN {S_{q^{1/2}}(\Sigma)^{(N)}}

Recall that
 $\fSiN$ is the field of fractions of $\SqSN$,   
$$\fSi=\SqS\otimes_{\SqSN}\fSiN,$$
 $K_{\Sigma,q^{1/2}}=\dim_{\fSiN}\fSi,$ and 
 $\lambda_{\Sigma,q^{1/2}}=\dim_{\SqSN} \SqS.$

%
%

\bl\label{upper}
Let $\Sigma$ be a pb surface,  we require $\chi(\Sigma)$ is negative if $\partial \Sigma = \emptyset$.
Then there exist $N^{3r(\Sigma)}$ elements in $S_{q^{1/2}}(\Sigma)$, which are linearly independent over 
$S_{q^{1/2}}(\Sigma)^{(N)}$. Especially we have $$N^{3r(\Sigma)}\leq K_{\Sigma,q^{1/2}}\leq \lambda_{\Sigma,q^{1/2}}.$$
\el
\bpr
 Lemma \ref{independence_skein}
indicates Lemma \ref{upper} when $\partial \Sigma= \emptyset$.
 When $\partial \Sigma\neq \emptyset$,
Lemmas \ref{Saturate} and \ref{independence} indicate Lemma \ref{upper}.
\epr

\bt\label{stated_u_l}
Let $\Sigma$ be a  pb surface with nonempty boundary.
Then we have 
\begin{align*}
N^{3 r(\Sigma)}\leq
\lambda_{\Sigma,q^{1/2}} \leq (2N^3 - \frac{N(N+1)(2N+1)}{6})^{r(\Sigma)}.
\end{align*}
\et
\bpr
Lemmas \ref{upper_O}, \ref{Saturate} and \ref{upper}.
\epr

\begin{corollary}
Let $\Sigma$ be a  pb surface with nonempty boundary.
Then we have
$$\lim_{N \text{ being odd }\rightarrow +\infty} \frac{\log \lambda_{\Sigma,q^{1/2}}}{\log N} = 3 r(\Sigma).$$
\end{corollary}

\bt\label{local_stated}
Let $\Sigma$ be a  pb surface, we require $\chi(\Sigma)$ is negative if $\partial \Sigma = \emptyset$.  
Then we have $K_{\Sigma,q^{1/2}}= N^{3 r(\Sigma)}$.
\et
\bpr
When $\partial = \emptyset$, Theorem \ref{local_stated} was proved in Theorem \ref{local_skein}.

Then assume $\partial\Sigma\neq\emptyset$.
From Lemma \ref{upper}, it suffices to show $K_{\Sigma,q^{1/2}}\leq N^{3 r(\Sigma)}.$
We use $\widetilde{(A_q)^{\otimes r(\Sigma)}}$ to denote the field of fractions of the commutative domain $(A_q)^{\otimes r(\Sigma)}$, 
and use
$ \widetilde{\Oq^{\otimes r(\Sigma)}} $ to denote $\Oq^{\otimes r(\Sigma)}\otimes_{(A_q)^{\otimes r(\Sigma)}} \widetilde{(A_q)^{\otimes r(\Sigma)}}$. 
Then from Lemma \ref{Saturate}, it suffices to show
$\dim_{\widetilde{(A_q)^{\otimes r(\Sigma)}}} \widetilde{\Oq^{\otimes r(\Sigma)}}\leq  N^{3 r(\Sigma)}.$
To simplify notation, we will regard elements in $\Oq^{\otimes r(\Sigma)}$
as elements in $ \widetilde{\Oq^{\otimes r(\Sigma)}} $ via the obvious embedding from 
 $\Oq^{\otimes r(\Sigma)}$
to $ \widetilde{\Oq^{\otimes r(\Sigma)}} $.

We use $x_1,\cdots,x_{N^3}$ to denote the basis elements in Theorem \ref{dimen}, then 
$x_1,\cdots,x_{N^3}$  are actually elements in $\Oq$. Let $V$ be the sub-vector space of $\widetilde{\Oq^{\otimes r(\Sigma)}}$ linearly spanned by 
$\{x_{i_1}\otimes\cdots\otimes x_{i_{r(\Sigma)}} \mid 1\leq i_1,\cdots,i_{r(\Sigma)}\leq N^3\}$.
Let $y_1,\cdots,y_{r(\Sigma)}\in \Oq$ be $r(\Sigma)$ elements in $\Oq$. For each $1\leq i\leq r(\Sigma)$, from Theorem \ref{dimen} there exist $u_i\in A_q-\{0\}$ and $v_{i,1},\cdots,v_{i,N^3}\in A_q$ such that 
$u_iy_i = \sum_{1\leq j\leq N^3} v_{i,j} x_j.$
Then we have 
\begin{align*}
&(u_1\otimes\cdots\otimes u_{r(\Sigma)})(y_1\otimes \cdots \otimes y_{r(\Sigma)})
= u_1y_1\otimes \cdots \otimes u_{r(\Sigma)}y_{r(\Sigma)}\\
=&\sum_{1\leq j_1,\cdots,j_{r(\Sigma)}\leq N^3} v_{1,j_1} x_{j_1}\otimes \cdots \otimes v_{r(\Sigma),j_{r(\Sigma)}}x_{r(\Sigma)}\\
=&\sum_{1\leq j_1,\cdots,j_{r(\Sigma)}\leq N^3}( v_{1,j_1} \otimes \cdots \otimes v_{r(\Sigma),j_{r(\Sigma)}})  (x_{j_1}\otimes \cdots \otimes x_{r(\Sigma)})
\end{align*}
where $u_1\otimes\cdots\otimes u_{r(\Sigma)}\in (A_q)^{\otimes r(\Sigma)}-\{0\}$ and 
$ v_{1,j_1} \otimes \cdots \otimes v_{r(\Sigma),j_{r(\Sigma)}}\in (A_q)^{\otimes r(\Sigma)}$.
Thus $y_1\otimes \cdots \otimes y_{r(\Sigma)}\in V$. Then clearly we have $V = \widetilde{\Oq^{\otimes r(\Sigma)}}$.
This completes the proof since the cardinality of $\{x_{i_1}\otimes\cdots\otimes x_{i_{r(\Sigma)}} \mid 1\leq i_1,\cdots,i_{r(\Sigma)}\leq N^3\}$ is $N^{3 r(\Sigma)}$.
\epr

\begin{remark}
When $\partial \Sigma\neq \emptyset$, obviously we have $$\{x_{i_1}\otimes\cdots\otimes x_{i_{r(\Sigma)}} \mid 1\leq i_1,\cdots,i_{r(\Sigma)}\leq N^3\}$$ in the  proof for Theorem \ref{local_stated} is actually a basis.
\end{remark}

\begin{remark}
For general stated $SL_n$-skein algebra, the Frobenius map was built when $\partial\Sigma\neq\emptyset$ \cite{wang2023stated}. From Theorem 7.3 in chapter 3 in \cite{brown2012lectures} for $SL_n$, and a general statement  for Lemma 6.28 for $SL_n$, using the same techniques we can generalize Theorem \ref{local_stated}  to $SL_n$, and the corresponding dimension is $N^{(n^2-1)r(\Sigma)}$.
\end{remark}

\subsection{When the marking set is nonempty}

For a compact marked three manifold $(M,\mathcal{N})$, we will give an upper bound for $$
\dim_{\SMN} \SMQ,$$  in this subsection. From Theorem \ref{reduced}, this offers an  upper bound for the dimension of the reduced stated skein module.

\bt
Let $(M,\mathcal{N})$ be a  compact  marked three manifold with $\mathcal{N}\neq \emptyset$.
Suppose the genus of $(M,\mathcal{N})$ is $g$, and $\mathcal{N}$ has $k$ components. Then we have 
$$\text{dim}_{\SMN}S_{q^{1/2}}(M,\mathcal{N})\leq (2N^3 - \frac{N(N+1)(2N+1)}{6})^{2g+k-1}.$$
\et
\bpr
Remark \ref{rem_upper}, Theorem \ref{stated_u_l}.
\epr

\bt
Let $(M,\mathcal{N})$ be a  compact  marked three manifold with $\mathcal{N}\neq \emptyset$.
Suppose the genus of $(M,\mathcal{N})$ is $g$, and $\mathcal{N}$ has $k$ components. Then we have the dimension of
any reduced stated skein module of $\MN$ is not more than
$$(2N^3 - \frac{N(N+1)(2N+1)}{6})^{2g+k-1}.$$
\et

\bibliographystyle{plain}

\bibliography{ref.bib}

\hspace*{\fill} \\

School of Physical and Mathematical Sciences, Nanyang Technological University, 21 Nanyang Link Singapore 637371

$\emph{Email address}$: zhihao003@e.ntu.edu.sg

\vspace{0.5cm}

Bernoulli Institute, University of Groningen, P.O. Box 407, 9700 AK Groningen, The Netherlands

$\emph{Email address}$: wang.zhihao@rug.nl

\end{document}

{\cred Try to calculate the dimension of $\Sq$ over $\SqN$ when $M$ is a handle body, more ambitious goal is when $M$ is the thickening of the surface (Techniques used in \cite{frohman2021dimension} maybe helpful). Try and study skein module over image of frobenius for lens space and periodic mapping Tori over torus, also when $M$ is the complement of two bridge knot or two bridge link.  $M$ is the thickening of the closed torus.

 Try to calculation the skein module of trivial  $S^1$-bundles over surfaces (the surface can have boundary) over the ring $\mathbb{C}[A^{\pm 1}]$ (maybe by using techniques in \cite{detcherry2021basis})

Maybe also try to  consider stated case. I think stated case is easier to handle with.

Theorems 2.2 and 2.3 in "Fundamentals of Kauffman bracket skein modules".

MULTIPLICATIVE STRUCTURE OF KAUFFMAN BRACKET
SKEIN MODULE QUANTIZATIONS

ALGEBRAIC GENERATORS OF THE SKEIN ALGEBRA OF A SURFACE

On the genus two skein algebra

On skein algebras of planar surfaces

Proposition 16 in REPRESENTATIONS OF THE KAUFFMAN SKEIN ALGEBRA OF SMALL SURFACES

The skein module of torus knots complements

SOME RESULTS ABOUT THE KAUFFMAN BRACKET SKEIN
MODULE OF THE TWIST KNOT EXTERIOR
}